\documentclass[12pt]{article}
\usepackage{amsmath,amsthm,amsfonts,amssymb,amscd}
\def\qed{{\unskip\nobreak\hfil\penalty50
\hskip2em\hbox{}\nobreak\hfil$\square$
\parfillskip=0pt \finalhyphendemerits=0\par}\medskip}
\def\proof{\trivlist \item[\hskip \labelsep{\bf Proof\ }]}
\def\endproof{\null\hfill\qed\endtrivlist}

\def\Ad{{\mathrm {Ad}}}

\def\Hom{{\mathrm {Hom}}}

\def\a{\alpha}
\def\b{\beta}

\def\e{\varepsilon}

\def\l{\langle}
\def\r{\rangle}

\def\phi{\varphi}

\def\Om{\Omega}

\def\emptyset{\varnothing}
\def\setminus{\smallsetminus}

\def\Diff{{\mathrm {Diff}}}
\def\Mob{{\rm\textsf{M\"ob}}}

\def\Ad{{\mathrm {Ad}}}

\def\Hom{{\mathrm {Hom}}}

\def\a{\alpha}
\def\b{\beta}

\def\e{\varepsilon}

\def\phi{\varphi}

\def\Om{\Omega}

\newtheorem{theorem}{Theorem}[section]
\newtheorem{lemma}[theorem]{Lemma}
\newtheorem{conjecture}[theorem]{Conjecture}
\newtheorem{corollary}[theorem]{Corollary}
\newtheorem{definition}[theorem]{Definition}

\newtheorem{proposition}[theorem]{Proposition}
\newtheorem{remark}[theorem]{Remark}

\def\emptyset{\varnothing}
\def\setminus{\smallsetminus}

\def\exp{{\mathrm {exp}}}

\def\Diff{{\mathrm {Diff}}}
\def\Mob{{\rm\textsf{M\"ob}}}

\def\res{\!\restriction\!}
\def\A{{\cal A}}

\def\B{{\cal B}}

\def\D{{\cal D}}

\def\I{{\cal I}}

\def\E{{\cal E}}
\def\H{{\cal H}}

\def\Z{{\mathbb Z}}

\renewcommand{\qed}{\ \hfill $\blacksquare$}

\newcommand{\bdefin}{\begin{definition}}
\newcommand{\blemma}{\begin{lemma}}
\newcommand{\bprop}{\begin{proposition}}
\newcommand{\btheor}{\begin{theorem}}
\newcommand{\bcoro}{\begin{corollary}}
\newcommand{\bconj}{\begin{conjecture}}
\newcommand{\edefin}{\end{definition}}
\newcommand{\elemma}{\end{lemma}}
\newcommand{\eprop}{\end{proposition}}
\newcommand{\etheor}{\end{theorem}}
\newcommand{\ecoro}{\end{corollary}}
\newcommand{\econj}{\end{conjecture}}
\newcommand{\brem}{\begin{remark}}
\newcommand{\erem}{\end{remark}}

\newcommand{\ba}{\begin{array}}
\newcommand{\ea}{\end{array}}

\renewcommand{\mod}{\mbox{mod}}

\title{\huge Some computations in the cyclic permutations of completely rational nets\\}
\author{
{\sc Feng Xu}\footnote{Supported in part by NSF.}\\
Department of Mathematics\\
University of California at Riverside\\
Riverside, CA 92521\\
E-mail: {\tt xufeng@math.ucr.edu}}
\begin{document}
\date{}
\maketitle

\begin{abstract}
In this paper we calculate certain chiral quantities from the cyclic
permutation orbifold of a general completely rational net. We
determine the fusion of a fundamental soliton, and by suitably
modified arguments of A. Coste , T. Gannon and especially P. Bantay
to our setting we are able to prove a number of arithmetic
properties including congruence subgroup properties for  $S, T$
matrices of a completely rational net defined by K.-H. Rehren .
2000MSC:81R15, 17B69.
\end{abstract}
\newpage

\section{Introduction}
Let $\A$ be a completely rational net (cf. Definition \ref{abr}).
Then $\A\otimes\A\otimes ...\otimes\A$ (n tensors) admits an action
of $\Bbb {Z}_n$ by cyclic permutations. The corresponding orbifold
net is referred to as the (n-th order) cyclic permutation orbifold
of $\A$. This construction has been used both in mathematics and
physics literature (for a partial list, see \cite {B1}, \cite{B2},
\cite{BDM}, \cite{BHS} and references therein). In \cite{LX}, this
construction was used for $n=2$ to show that strong additivity is
automatic in a conformal net with finite $\mu$ index. In \cite{KLX},
the construction was used to demonstrate applications of general
orbifold theories among other things. The starting motivation of
this paper is to improve a result (Prop. 9.4 in \cite{KLX}) on
fusion of a fundamental soliton. The second motivation came from two
papers: one by A. Coste and T. Gannon (cf. \cite{CG}) where under
certain conditions they showed that the $S,T$ matrices verified
congruence subgroup properties, and one by P. Bantay (cf.
\cite{B1})where he showed that congruence subgroup properties hold
if a number a heuristic arguments including what he called
``Orbifold Covariance Principle" hold. This ``Orbifold Covariance
Principle" of P. Bantay is highly nontrivial even in concrete
examples, and at present the only conceptual framework in which this
principle is a theorem is in the framework of local conformal net
(cf. Section \ref{nets}), where the principle follows from Theorem
\ref{orb}. In the language of local conformal nets the $S,T$
matrices were defined by K.-H. Rehren (cf. \cite{R2}) by using local
data of conformal nets, and in all known cases they agree with the
``S,T" matrices coming from modular transformations of characters.
If one is interested in the modular tensor categories, then this
``S,T" matrices of Rehren are sufficient for calculations of three
manifold invariants (cf. for example \cite{X2}). It is therefore an
interesting question to see if one can adapt the methods of A. Coste
and T. Gannon and P. Bantay to Rehren's ``S,T" matrices. \par In
this paper we will show that Prop. 9.4 of \cite{KLX} holds in
general (cf. Th. \ref{94}), and that a suitable modification of the
arguments of A. Coste and T. Gannon and P. Bantay is possible in the
conformal net setting, and congruence subgroup properties hold for
Rehern's ``S,T" matrices (cf. Th. \ref{cong}). Our key observation
is the  squares of nets in \S3. By using (3) of Lemma \ref{Smatrix1}
which relates the chiral data of a net and subnet for suitably
chosen squares, we are able to obtain strong constraints on certain
matrices (cf. Th. \ref{keyeq}). These squares are in fact commuting
squares first considered in the setting of $II_1$ factors by S. Popa
in \cite{Po}, and they already played an important role in the
setting of nets in \cite{X5}. However the ``commuting" property of
these squares will not play an explicit role in this paper.  Th.
\ref{keyeq} allows us to apply the methods of P. Bantay in \cite{B1}
to obtain arithmetic properties of Rehren's ``S,T" matrices in Th.
\ref{arithmetic1} and Th. \ref{cong}.

We note that Th. \ref{keyeq} implies series of arithmetic properties
of ``S,T" matrices, and even the first one as observed in \cite{KLX}
seems to be nontrivial  for concrete examples like the nets coming
from $SU(n)$ at level $k$.
\par The rest of this paper is organized as follows: In \S2,
 after recalling basic definitions of completely rational nets,
Rehren's $S,T$ matrices, orbifolds and Galois actions, we stated a
few general results from \S9 of \cite{KLX} which will be used in \S3
and \S4. In \S3 we improve Prop. 9.4 of \cite{KLX} in Th. \ref{94},
and  we present the proof of Th. \ref{keyeq} as mentioned above from
a commuting square. In \S4, by modifications incorporating phase
factors the arguments of A. Coste and T. Gannon and P. Bantay, we
are able to prove Th. \ref{arithmetic1} and Th. \ref{cong}. We note
that the arguments in \S4 can be simplified if one can prove a
conjecture on Page 734 of \cite{KLX}.

\par

\section{Conformal nets, complete rationality, and orbifolds}

For the convenience of the reader we collect here some basic notions
that appear in this paper. This is only a guideline and the reader
should look at the references for a more complete treatment.

\subsection{Conformal nets on $S^1$} \label{nets}

By an interval of the circle we mean an open connected
non-empty subset $I$ of $S^1$ such that the interior of its
complement $I'$ is not empty.
We denote by $\I$ the family of all intervals of $S^1$.

A {\it net} $\A$ of von Neumann algebras on $S^1$ is a map
\[
I\in\I\to\A(I)\subset B(\H)
\]
from $\I$ to von Neumann algebras on a fixed Hilbert space $\H$
that satisfies:
\begin{itemize}
\item[{\bf A.}] {\it Isotony}. If $I_{1}\subset I_{2}$ belong to
$\I$, then
\begin{equation*}
 \A(I_{1})\subset\A(I_{2}).
\end{equation*}
\end{itemize}
The net $\A$ is called {\it local} if it satisfies:
\begin{itemize}
\item[{\bf B.}] {\it Locality}. If $I_{1},I_{2}\in\I$ and $I_1\cap
I_2=\emptyset$ then
\begin{equation*}
 [\A(I_{1}),\A(I_{2})]=\{0\},
 \end{equation*}
where brackets denote the commutator.
\end{itemize}
The net $\A$ is called {\it M\"{o}bius covariant} if in addition
satisfies
the following properties {\bf C,D,E,F}:
\begin{itemize}
\item[{\bf C.}] {\it M\"{o}bius covariance}.
There exists a strongly
continuous unitary representation $U$ of the M\"{o}bius group
$\Mob$ (isomorphic to $PSU(1,1)$) on $\H$ such that
\begin{equation*}
 U(g)\A(I) U(g)^*\ =\ \A(gI),\quad g\in \Mob,\ I\in\I.
\end{equation*}
\end{itemize}
If $E\subset S^1$ is any region, we shall put
$\A(E)\equiv\bigvee_{E\supset I\in\I}\A(I)$ with $\A(E)=\mathbb C$ if
$E$ has empty interior (the symbol $\vee$ denotes the von Neumann
algebra generated).  Note that the definition of $\A(E)$ remains the
same if $E$ is an interval namely: if $\{I_n\}$ is an increasing
sequence of intervals and $\cup_n I_n = I$, then the $\A(I_n)$'s generate
$\A(I)$ (consider a sequence of elements $g_n\in\Mob$ converging to the
identity such that $g_n I\subset I_n$).
\begin{itemize}
\item[{\bf D.}] {\it Positivity of the energy}.
The generator of the one-parameter
rotation subgroup of $U$ (conformal Hamiltonian) is positive.
\item[{\bf E.}] {\it Existence of the vacuum}.  There exists a unit
$U$-invariant vector $\Omega\in\H$ (vacuum vector), and $\Omega$ is
cyclic for the von Neumann algebra $\bigvee_{I\in\I}\A(I)$.
\end{itemize}
By the Reeh-Schlieder theorem $\Omega$ is cyclic and separating for
every fixed $\A(I)$. The modular objects associated with
$(\A(I),\Omega)$ have a geometric meaning
\[
\Delta^{it}_I = U(\Lambda_I(2\pi t)),\qquad J_I = U(r_I)\ .
\]
Here $\Lambda_I$ is a canonical one-parameter subgroup of $\Mob$ and $U(r_I)$ is a
antiunitary acting geometrically on $\A$ as a reflection $r_I$ on $S^1$.

This implies {\em Haag duality}:
\[
\A(I)'=\A(I'),\quad I\in\I\ ,
\]
where $I'$ is the interior of $S^1\setminus I$.

\begin{itemize}
\item[{\bf F.}] {\it Irreducibility}. $\bigvee_{I\in\I}\A(I)=B(\H)$.
Indeed $\A$ is irreducible iff
$\Om$ is the unique $U$-invariant vector (up to scalar multiples).
Also  $\A$ is irreducible
iff the local von Neumann
algebras $\A(I)$ are factors. In this case they are III$_1$-factors in
Connes classification of type III factors
(unless $\A(I)=\mathbb C$ for all $I$).
\end{itemize}
By a {\it conformal net} (or diffeomorphism covariant net)
$\A$ we shall mean a M\"{o}bius covariant net such that the following
holds:
\begin{itemize}
\item[{\bf G.}] {\it Conformal covariance}. There exists a projective
unitary representation $U$ of $\Diff(S^1)$ on $\H$ extending the unitary
representation of $\Mob$ such that for all $I\in\I$ we have
\begin{gather*}
 U(g)\A(I) U(g)^*\ =\ \A(gI),\quad  g\in\Diff(S^1), \\
 U(g)xU(g)^*\ =\ x,\quad x\in\A(I),\ g\in\Diff(I'),
\end{gather*}
\end{itemize}
where $\Diff(S^1)$ denotes the group of smooth, positively oriented
diffeomorphism of $S^1$ and $\Diff(I)$ the subgroup of
diffeomorphisms $g$ such that $g(z)=z$ for all $z\in I'$.
\par
Let $G$ be a simply connected  compact Lie group. By Th. 3.2 of
\cite{FG}, the vacuum positive energy representation of the loop
group $LG$ (cf. \cite{PS}) at level $k$ gives rise to an irreducible
conformal net denoted by {\it ${\A}_{G_k}$}. By Th. 3.3 of
\cite{FG}, every irreducible positive energy representation of the
loop group $LG$ at level $k$ gives rise to  an irreducible covariant
representation of ${\A}_{G_k}$.

\subsection{Doplicher-Haag-Roberts superselection sectors in CQFT}
The references of this section are \cite{DHR,FRS,L1,L1',LR, GL1,
GL2}. The DHR theory was originally made on the 4-dimensional
Minkowski spacetime, but can be generalized to our setting. There
are however several important structure differences in the low
dimensional case.

A (DHR) representation $\pi$ of $\A$ on a Hilbert space $\H$ is a map
$I\in\I\mapsto  \pi_I$ that associates to each $I$ a normal
representation of $\A(I)$ on $B(\H)$ such that
\[
\pi_{\tilde I}\res\A(I)=\pi_I,\quad I\subset\tilde I, \quad
I,\tilde I\subset\I\ .
\]
$\pi$ is said to be M\"obius (resp. diffeomorphism) covariant if
there is a projective unitary representation $U_{\pi}$ of $\Mob$ (resp.
$\Diff^{(\infty)}(S^1)$, the infinite cover of $\Diff(S^1)$ ) on $\H$ such that
\[
\pi_{gI}(U(g)xU(g)^*) =U_{\pi}(g)\pi_{I}(x)U_{\pi}(g)^*
\]
for all $I\in\I$, $x\in\A(I)$ and $g\in \Mob$ (resp.
$g\in\Diff^{(\infty)}(S^1)$). Note that if $\pi$ is irreducible and
diffeomorphism covariant then $U$ is indeed a projective unitary
representation of $\Diff(S^1)$.

By definition the irreducible conformal net is in fact an irreducible
representation of itself and we will call this representation the {\it
vacuum representation}.

Given an interval $I$ and a representation $\pi$ of $\A$, there is
an {\em endomorphism of $\A$ localized in $I$} equivalent to $\pi$;
namely $\rho$ is a representation of $\A$ on the vacuum Hilbert
space $\H$, unitarily equivalent to $\pi$, such that
$\rho_{I'}=\text{id}\restriction\A(I')$.

Fix an interval $I_0$ and  endomorphisms $\rho,\rho'$ of $\A$
localized in $I_0$. Then the {\em composition} (tensor product)
$\rho\rho'$ is defined by
\[
(\rho\rho')_I=\rho_I \rho'_I
\]
with $I$ an interval containing $I$. One can indeed define
$(\rho\rho')_I$ for an arbitrary interval $I$ of $S^1$ (by using
covariance) and get a well defined endomorphism of $\A$ localized in
$I_0$. Indeed the endomorphisms of $\A$ localized in a given
interval form a tensor $C^*$-category. For our needs $\rho,\rho'$
will be always localized in a common interval $I$.

If $\pi$ and $\pi'$ are representations of $\A$, fix an interval
$I_0$ and choose endomorphisms $\rho,\rho'$ localized in $I_0$ with
$\rho$ equivalent to $\pi$ and $\rho'$ equivalent to $\pi'$. Then
$\pi\cdot\pi'$ is defined (up to unitary equivalence)  to be
$\rho\rho'$. The class of a DHR representation modulo unitary
equivalence is a {\em superselection sectors} (or simply a sector).
We use the notations $\rho_1\succ \rho_2$ or $\rho_2\prec \rho_1$ if
$\rho_2$ appears a summand of $\rho_1.$

The localized endomorphisms of $\A$ for a tensor $C^*$-category. For
our needs, $\rho,\rho'$ will be always localized in a common
interval $I$.

We now define  the statistics. Given the endomorphism $\r$ of $\A$
localized in $I\in\I$, choose an equivalent endomorphism $\r_0$
localized in an interval $I_0\in\I$ with $\bar I_0\cap\bar I
=\emptyset$ and let $u$ be a local intertwiner in $\Hom(\r,\r_0)$ as
above, namely $u\in \Hom(\rho_{\tilde I},\rho_{0,\tilde I})$ with
$I_0$ following clockwise $I$ inside $\tilde I$ which is an interval
containing both $I$ and $I_0$.

The {\it statistics operator} $\varepsilon := u^*\rho(u) =
u^*\rho_{\tilde I}(u) $ belongs to $\Hom(\rho^2_{\tilde
I},\rho^2_{\tilde I})$.  An elementary computation shows that it
gives rise to a presentation of  the Artin  braid group
$$\epsilon_i\epsilon_{i+1}\epsilon_i =
\epsilon_{i+1}\epsilon_i\epsilon_{i+1},\qquad
\epsilon_i\epsilon_{i'} =\epsilon_{i'}\epsilon_i \,\quad{\rm if}\,\,
|i-i'|\geq 2,$$ where $\varepsilon_i=\rho^{i-1}(\varepsilon)$. The
(unitary equivalence class of the) representation of the Artin braid
group thus obtained is the {\it statistics} of the superselection
sector $\rho$.

It turns out the endomorphisms localized in a given interval form a {\em
braided $C^*$-tensor category} with unitary braiding.

The {\em statistics parameter} $\lambda_\rho$ can be defined in
general. In particular, assume $\rho$ to be localized in $I$ and
$\rho_I\in\text{End}((\A(I))$ to be irreducible with a conditional
expectation $E: \A(I)\to \rho_I(\A(I))$, then
\[
\lambda_\rho:=E(\epsilon)
\]
depends only on the superselection sector of $\rho$.

The {\em statistical dimension} $d(\rho)$ and the  {\it univalence}
$\omega_\rho$ are then defined by
\[
d(\rho) = |\lambda_\rho|^{-1}\ ,\qquad \omega_\rho =
\frac{\lambda_\rho}{|\lambda_\rho|}\ .
\]
The {\em conformal spin-statistics theorem} shows that
\[
\omega_\rho = e^{i 2\pi \Delta_\rho}\ ,
\]
where $\Delta_\rho$ is the conformal dimension (the lowest
eigenvalue of the generator of the rotation subgroup) in the
representation $\rho$. The right hand side in the above equality is
called the {\em univalence} of $\rho$.\par $d(\rho)^2$ will be
called the index of $\rho$. The general index  was first defined and
investigated  by Vaughan Jones in the case of $II_1$ factors in
\cite{J}.


\subsection{Rehren's  $S,T$-matrices}
Next we will recall some of the results of \cite{R2}  and introduce
notations. \par Let $\{[\lambda], \lambda \in \mathcal{P} \}$ be a
finite set of all equivalence classes of irreducible, covariant,
finite-index representations of an irreducible local conformal net
$\A$. We will denote the conjugate of $[\lambda]$ by $[{\bar
\lambda}]$ and identity sector (corresponding to the vacuum
representation) by $[1]$ if no confusion arises, and let
$N_{\lambda\mu}^\nu = \langle [\lambda][\mu], [\nu]\rangle $. Here
$\langle \mu,\nu\rangle$ denotes the dimension of the space of
intertwiners from $\mu$ to $\nu$ (denoted by $\text {\rm
Hom}(\mu,\nu)$).  We will denote by $\{T_e\}$ a basis of isometries
in $\text {\rm Hom}(\nu,\lambda\mu)$. The univalence of $\lambda$
and the statistical dimension of (cf. \S2  of \cite{GL1}) will be
denoted by $\omega_{\lambda}$ and $d{(\lambda)}$ (or $d_{\lambda})$)
respectively. \par Let $\phi_\lambda$ be the unique minimal left
inverse of $\lambda$, define:
\begin{equation}\label{Ymatrix}
Y_{\lambda,\mu}:= d(\lambda)  d(\mu) \phi_\mu (\epsilon (\mu,
\lambda)^* \epsilon (\lambda, \mu)^*),
\end{equation}
where $\epsilon (\mu, \lambda)$ is the unitary braiding operator
 (cf. \cite{GL1} ). \par
We list two properties of $Y_{\lambda, \mu}$ (cf. (5.13), (5.14) of
\cite{R2}) which will be used in the following:
\begin{lemma}\label{Yprop}
\begin{equation*}
Y_{\lambda,\mu} = Y_{\mu,\lambda}  = Y_{\lambda,\bar \mu}^* =
Y_{\bar \lambda, \bar \mu}.
\end{equation*}
\begin{equation*}
Y_{\lambda,\mu}  = \sum_k N_{\lambda\mu}^\nu
\frac{\omega_\lambda\omega_\mu} {\omega_\nu} d(\nu) .
\end{equation*}
\end{lemma}
We note that one may take the second equation in the above lemma as
the definition of $Y_{\lambda,\mu}$.\par Define $a := \sum_i
d_{\rho_i}^2 \omega_{\rho_i}^{-1}$. If the matrix $(Y_{\mu,\nu})$ is
invertible, by Proposition on P.351 of \cite{R2} $a$ satisfies
$|a|^2 = \sum_\lambda d(\lambda)^2$.
\begin{definition}\label{c0}
Let $a= |a| \exp(-2\pi i \frac{c_0(\A)}{8})$ where  $c_0(\A)\in
{\mathbb R}$ and $c_0(\A)$ is well defined ${\rm mod} \ 8\mathbb Z$.
For simplicity we will denote $c_0(\A)$ simply as $c_0$ when the
underlying $\A$ is clear.
\end{definition}
Define matrices
\begin{equation}\label{Smatrix}
S:= |a|^{-1} Y, T:=  C {\rm Diag}(\omega_{\lambda})
\end{equation}
where \[C:= \label{dims} \exp(-2\pi i \frac{c_0}{24}).\] Then these
matrices satisfy (cf. \cite{R2}):
\begin{lemma}\label{Sprop}
\begin{align*}
SS^{\dag} & = TT^{\dag} ={\rm id},  \\
STS &= T^{-1}ST^{-1},  \\
S^2 & =\hat{C},\\
 T\hat{C} & =\hat{C}T=T,
\end{align*}

where $\hat{C}_{\lambda\mu} = \delta_{\lambda\bar \mu}$
is the conjugation matrix.
\end{lemma}

The above Lemma shows that $S,T$ as defined there give rise to a
representation of the modular group denoted by $\Gamma(1).$ This is
the group generated by two matrices \( t=\left( \begin{array}{cc}
1 & 1\\
0 & 1
\end{array}\right)  \) and \( s=\left( \begin{array}{cc}
0 & -1\\
1 & 0
\end{array}\right)  \), and the representation is given by
$s\rightarrow S, t\rightarrow T.$

Let $r$ be a rational number. Throughout this paper we will use
$T^r$ to denote a diagonal matrix whose $(\lambda,\lambda)$ entry is
given by $ \exp(2\pi i(\Delta_\lambda-\frac{c_0}{24})r).$

Moreover
\begin{equation}\label{Verlinde}
N_{\lambda\mu}^\nu = \sum_\delta \frac{S_{\lambda,\delta}
S_{\mu,\delta} S_{\nu,\delta}^*}{S_{1,\delta}}. \
\end{equation}
is known as Verlinde formula. \par Sometimes we will refer the $S,T$
matrices as defined above  as  {\bf genus 0 modular matrices of
${\A}$} since they are constructed from the fusion rules,
monodromies and minimal indices which can be thought as  genus 0
{\bf chiral data} associated to a Conformal Field Theory. \par Let
$c$ be the central charge associated with the projective
representations of ${\rm Diff}(S^1)$ of the conformal net $\A$ (cf.
\cite{K1} ). Note that by \cite{CW} $c$ is uniquely determined for a
conformal net. We will see that $c$ is always rational for a
completely rational net (see (4) of Th. \ref{arithmetic1} for a more
refined statement).

It is   proved in Lemma 9.7 of \cite{KLX} that $c_0-c\in 4\mathbb Z$
for complete rational nets. \par The commutative algebra generated
by $\lambda$'s with structure constants $N_{\lambda\mu}^\nu$ is
called {\bf fusion algebra} of $\A$. If $Y$ is invertible, it
follows from Lemma \ref{Sprop} and equation (\ref{Verlinde}) that
any nontrivial irreducible representation of the fusion algebra is
of the form $\lambda\rightarrow \frac{S_{\lambda\mu}}{S_{1\mu}}$ for
some $\mu$.
\par
\subsection{The Galois action on Rehren's $S,T$ matrices}
The basic idea in the theory of the Galois action \cite{BG}\cite{CG}
is to look at the field \( F \) obtained by adjoining to the
rationales \( \mathbb {Q} \) the matrix elements of all modular
transformations as defined after Lemma \ref{Sprop}. One may show
that, as a consequence of Verlinde's formula, \( F \) is a finite
Abelian extension of \( \mathbb {Q} \). By the  theorem of Kronecker
and Weber this means that \( F \) is a subfield of some cyclotomic
field \( \mathbb {Q}\left[ \zeta _{n}\right]  \) for some integer \(
n \), where \( \zeta _{n}=\exp \left( \frac{2\pi i}{n}\right)  \) is
a primitive \( n \)-th root of unity. We'll call the conductor of \(
\mathcal{A} \) the smallest \( n \) for which \( F\subseteq \mathbb
{Q}\left[ \zeta _{n}\right]  \) and which is divisible by the order
of the $T$ matrix.

The above results imply that the Galois group \( \mathrm{Gal}\left(
F/\mathbb {Q}\right)  \) is a homomorphic image of the Galois group
\( \mathcal{G}_{n}=\mathrm{Gal}\left(  \mathbb {Q}\left[ \zeta
_{n}\right] /\mathbb {Q}\right)  \). But it is  known that \(
\mathcal{G}_{n} \) is isomorphic to the group \( \left( \Z /n\Z
\right) ^{*} \) of prime residues modulo \( n \), its elements being
the Frobenius maps \( \sigma _{l}:\mathbb {Q}\left[ \zeta
_{n}\right]\rightarrow \left[ \zeta _{n}\right] \) that leave \(
\mathbb {Q} \) fixed, and send \( \zeta _{n} \) to \( \zeta _{n}^{l}
\) for \( l \) coprime to \( n \). Consequently, the maps \( \sigma
_{l} \) are automorphisms of \( F \) over \( \mathbb {Q} \).

According to \cite{CG}, we have (for \( l \) coprime to the
conductor) \begin{equation} \label{sls} \sigma _{l}\left(
S_{\lambda,\mu}\right) =\varepsilon _{l}(\mu)S_{\lambda,\pi
_{l}(\mu)}
\end{equation}
for some permutation \( \pi _{l}\in Sym\left( \mathcal{P}\right)  \)
of the irreducible representations and some function \( \varepsilon
_{l}:\mathcal{P}\rightarrow \left\{ -1,+1\right\}  \). Upon
introducing the orthogonal monomial matrices
\begin{equation} \label{gldef} \left( G_{l}\right)
_{\lambda,\mu}=\varepsilon _{l}(\mu)\delta _{\lambda,\pi _{l}(\mu)}
\end{equation}
 and denoting by \( \sigma _{l}\left( M\right)  \) the matrix that one obtains
by applying \( \sigma _{l} \) to \( M \) elementwise, we have
\begin{equation} \label{slg} \sigma _{l}\left( S\right)
=SG_{l}=G_{l}^{-1}S
\end{equation}
Note that for \( l \) and \( m \) both coprime to the conductor
\begin{eqnarray*}
\pi _{lm} & = & \pi _{l}\pi _{m}\\
G_{lm} & = & G_{l}G_{m}
\end{eqnarray*}
 The Galois action on \( T \) is given by
 \begin{equation} \label{slt} \sigma _{l}\left( T\right) =T^{l}
\end{equation}

\subsection{The orbifolds}
Let ${\A}$ be an irreducible conformal net on a Hilbert space
${\H}$ and let $\Gamma$ be a finite group. Let $V:\Gamma\rightarrow U({\H})$
be a  unitary representation of $\Gamma$ on ${\H}$.
If $V:\Gamma\rightarrow U({\H})$ is not faithful, we set $\Gamma':= \Gamma
/{\rm ker} V$.
\begin{definition} \label{p'}
We say that $\Gamma$ acts properly on ${\A}$ if the following conditions
are satisfied:\par
(1) For each fixed interval $I$ and each $g\in \Gamma$,
$\alpha_g (a):=V(g)aV(g^*) \in {\A}(I), \forall a\in
{\A}(I)$; \par
(2) For each  $g\in \Gamma$, $V(g)\Omega = \Omega, \forall g\in \Gamma$.\par
\label{Definition 2.1}
\end{definition}
We note that if $\Gamma$ acts properly, then $V(g)$, $g\in\Gamma$
commutes with the unitary representation $U$ of $\Mob$. \par Define
$\B(I):= \{ a\in \A(I) | \alpha_g (a)=a, \forall g\in \Gamma \}$ and
${\A}^\Gamma(I):={\B}(I)P_0$ on ${\H}_0$ where $\H_0:=\{ x\in \H|
V(g)x=x, \forall g\in \Gamma \}$ and $P_0$ is the projection from
$\H$ to $\H_0.$  Then $U$ restricts to an  unitary representation
(still denoted by $U$) of $\Mob$ on ${\H}_0$. Then: \bprop The map
$I\in {\I}\rightarrow {\A}^{\Gamma}(I)$ on $ {\H}_0$ together with
the  unitary representation (still denoted by $U$) of \Mob\ on
${\H}_0$ is an irreducible M\"{o}bius covariant net.
\label{Prop.2.1} \eprop The irreducible  M\"{o}bius covariant net in
Prop.  \ref{Prop.2.1} will be denoted by ${\A}^\Gamma$ and will be
called the {\it orbifold of ${\A}$} with respect to $\Gamma$. When
$\Gamma$ is generated by $h_1,...,h_k$, we will write ${\A}^\Gamma$
as $\A^{\l h_1,...,h_k\r}.$

\par
\subsection{Complete rationality }
\label{complete rationality} We first recall some definitions from
\cite{KLM} . Recall that   ${\I}$ denotes the set of intervals of
$S^1$. Let $I_1, I_2\in {\I}$. We say that $I_1, I_2$ are disjoint
if $\bar I_1\cap \bar I_2=\emptyset$, where $\bar I$ is the closure
of $I$ in $S^1$. Denote by ${\I}_2$ the set of unions of disjoint 2
elements in ${\I}$. Let ${\A}$ be an irreducible M\"{o}bius
covariant net as in \S2.1. For $E=I_1\cup I_2\in{\I}_2$, let
$I_3\cup I_4$ be the interior of the complement of $I_1\cup I_2$ in
$S^1$ where $I_3, I_4$ are disjoint intervals. Let
$$
{\A}(E):= A(I_1)\vee A(I_2), \quad
\hat {\A}(E):= (A(I_3)\vee A(I_4))'.
$$ Note that ${\A}(E) \subset \hat {\A}(E)$.
Recall that a net ${\A}$ is {\it split} if ${\A}(I_1)\vee {\A}(I_2)$
is naturally isomorphic to the tensor product of von Neumann
algebras ${\A}(I_1)\otimes {\A}(I_2)$ for any disjoint intervals
$I_1, I_2\in {\I}$. ${\A}$ is {\it strongly additive} if
${\A}(I_1)\vee {\A}(I_2)= {\A}(I)$ where $I_1\cup I_2$ is obtained
by removing an interior point from $I$. \bdefin\label{abr}
\cite{KLM} ${\A}$ is said to be completely  rational if ${\A}$ is
split, strongly additive, and the index $[\hat {\A}(E): {\A}(E)]$ is
finite for some $E\in {\I}_2$ . The value of the index $[\hat
{\A}(E): {\A}(E)]$ (it is independent of $E$ by Prop. 5 of
\cite{KLM}) is denoted by $\mu_{{\A}}$ and is called the $\mu$-index
of ${\A}$. If the index $[\hat {\A}(E): {\A}(E)]$ is infinity for
some $E\in {\I}_2$, we define the $\mu$-index of ${\A}$ to be
infinity. \label{Definition 2.2} \edefin A formula for the
$\mu$-index of a subnet is proved in \cite{KLM}. With the result on
strong additivity for $\A^{\Gamma}$ in \cite{X1}, we have the
complete rationality in following theorem.

Note that, by our recent results in \cite{LX}, every irreducible,
split, local conformal net with finite $\mu$-index is automatically
strongly additive. \btheor\label{orb} Let ${\A}$ be an irreducible
M\"{o}bius covariant net and let $\Gamma$ be a finite group acting
properly on ${\A}$. Suppose that ${\A}$ is completely rational.
Then:\par (1): ${\A}^\Gamma$ is completely rational  and
$\mu_{{\A}^\Gamma}= |\Gamma'|^2 \mu_{{\A}}$; \par (2): There are
only a finite number of irreducible covariant representations of
${\A}^\Gamma$ (up to unitary equivalence), and they give rise to a
unitary modular category as defined in II.5 of \cite{Tu} by the
construction as given in \S1.7 of \cite{X2}. \label{Th.2.6} \etheor
Suppose that ${\A}$ and $\Gamma$ satisfy the assumptions of Th.
\ref{orb}. Then ${\A}^\Gamma$ has only finitely number of
irreducible representations $\dot\lambda$ and
$$
\sum_{\dot\lambda}d(\dot\lambda)^2 = \mu_{{\A}^\Gamma}= |\Gamma'|^2 \mu_{{\A}} .
$$

The set of such $\dot\lambda$'s is closed under conjugation and
compositions, and by Cor. 32 of \cite{KLM}, the $Y$-matrix in
(\ref{Ymatrix}) for ${\A}^\Gamma$ is non-degenerate, and we will
denote the corresponding genus $0$ modular matrices by $\dot S, \dot
T$. Denote by $\dot\lambda$ (resp. $\mu$) the irreducible covariant
representations of ${\A}^\Gamma$ (resp. ${\A}$) with finite index.
Denote by $b_{\mu\dot\lambda}\in {\mathbb N}\cup\{0\}$ the
multiplicity of representation $\dot\lambda$ which appears in the
restriction of representation $\mu$ when restricting from $ {\A}$ to
${\A}^{\Gamma}$. The $b_{\mu\dot\lambda} $ are also known as the
{\em branching rules.} An irreducible covariant representation
$\dot\lambda$  of ${\A}^\Gamma$ is called an {\it untwisted}
representation if $b_{\mu\dot\lambda}\neq 0$ for some representation
$\mu$ of $\A$. These are representations of ${\A}^\Gamma$ which
appear as subrepresentations in the restriction of some
representation of ${\A}$ to ${\A}^\Gamma$. A representation is
called {\it twisted} if it is not untwisted.
\subsection{Restriction to the real line: Solitons}
Denote by $\I_0$ the set of open, connected, non-empty, proper
subsets of $\mathbb R$, thus $I\in\I_0$ iff $I$ is an open interval
or half-line (by an interval of $\mathbb R$ we shall always mean a
non-empty open bounded interval of $\mathbb R$).

Given a net $\A$ on $S^1$ we shall denote by $\A_0$ its restriction
to $\mathbb R = S^1\setminus\{-1\}$. Thus $\A_0$ is an isotone map
on $\I_0$, that we call a \emph{net on $\mathbb R$}. In this paper
we denote by $J_0:=(0,\infty)\subset \mathbb R$.

A representation $\pi$ of $\A_0$ on a Hilbert space $\H$ is a map
$I\in\I_0\mapsto\pi_I$ that associates to each $I\in\I_0$ a normal
representation of $\A(I)$ on $B(\H)$ such that
\[
\pi_{\tilde I}\res\A(I)=\pi_I,\quad I\subset\tilde I, \quad I,\tilde
I\in\I_0\ .
\]
A representation $\pi$ of $\A_0$ is also called a \emph{soliton}. As
$\A_0$ satisfies half-line duality, namely
$$
\A_0(-\infty,a)'= \A_0(a,\infty), a\in \mathbb R,
$$
by the usual DHR argument \cite{DHR} $\pi$ is unitarily equivalent
to a representation $\rho$ which  acts identically on
$\A_0(-\infty,0)$, thus $\rho$ restricts to an endomorphism of
$\A(J_0)= \A_0(0,\infty)$. $\rho$ is said to be localized on $J_0$
and we also refer to $\rho$ as soliton endomorphism.

Clearly a representation $\pi$ of $\A$ restricts to a soliton
$\pi_0$ of $\A_0$. But a representation $\pi_0$ of $\A_0$ does not
necessarily extend to a representation of $\A$.

If $\A$ is strongly additive, and a representation $\pi_0$ of $\A_0$
extends to a DHR representation of $\A$, then it is easy to see that
such an extension is unique, and in this case we will use the same
notation $\pi_0$ to denote the corresponding DHR  representation of
$\A$.
\subsection{Induction and restriction for a net and its subnet}

Let $\A$ be a M\"obius covariant net. By a M\"obius (resp.
conformal) covariant subnet $\B\subset \A$ we mean a map
$$I\in \I\rightarrow \B(I)\subset\A(I)$$
that associates to each $I\in \I$ a von Neumann subalgebra $\B(I)$
so that isotony and covariance with respect to the M\"obius  (resp.
conformal) group  hold.

Given  a bounded interval $I_0\in\I_0$ we fix canonical endomorphism
$\gamma_{I_0}$ associated with $\B(I_0)\subset\A(I_0)$. Then we can
choose for each $I\subset\I_0$ with $I\supset I_0$ a canonical
endomorphism $\gamma_{I}$ of $\A(I)$ into $\B(I)$ in such a way that
$\gamma_{I}\res\A(I_0)=\gamma_{I_0}$ and $\gamma_{I_1}$ is the
identity on $\B(I_1)$ if $I_1\in\I_0$ is disjoint from $I_0$, where
$\gamma_{I}\equiv\gamma_{I}\res\B(I)$.

We then have an endomorphism $\gamma$ of the $C^*$-algebra
$\mathfrak A\equiv\overline{\cup_{I}\A(I)}$
($I$ bounded interval of $\mathbb R$).

Given a DHR endomorphism $\rho$ of $\B$ localized in $I_0$, the
induction $\a_{\rho}$ of $\rho$ is the endomorphism of $\mathfrak A$
given by
\[
\a_{\rho}\equiv
\gamma^{-1}\cdot\Ad\e(\rho,\gamma)\cdot\rho\cdot\gamma\ ,
\]
where $\e$ denotes the right braiding unitary symmetry (there is
another choice for $\a$ associated with the left braiding).
$\a_{\rho}$ is localized in a right half-line containing $I_0$,
namely $\a_\rho$ is the identity on $\A(I)$ if $I$ is a bounded
interval contained in the left complement of $I_0$ in $\mathbb R$.
Up to unitarily equivalence, $\a_\rho$ is localizable in any right
half-line thus $\a_\rho$ is normal on left half-lines, that is to
say, for every $a\in\mathbb R$, $\a_\rho$ is normal on the
$C^*$-algebra $\mathfrak A(-\infty,a)\equiv\overline{\cup_{I\subset
(-\infty,a)}\A(I)}$ ($I$ bounded interval of $\mathbb R$), namely
$\a_\rho\res\mathfrak A(-\infty,a)$ extends to a normal morphism of
$\A(-\infty,a)$. When there are several subnets involved, we will
use notation $\a^{\B\rightarrow \A}_\rho$  introduced in \S3 of
\cite{X4} to indicate the net and subnet where we apply the
induction.

\subsection{Preliminaries on cyclic orbifolds}\label{cyc}
In the rest of this paper we assume that $\A$ is completely
rational. $\D:= \A\otimes\A...\otimes\A$ ($n$-fold tensor product)
and $\B:=\D^{\mathbb Z_n}$ (resp. $\D^{\mathbb P_n}$ where $\mathbb
P_n$ is the permutation group on $n$ letters) is the fixed point
subnet of $\D$ under the action of cyclic permutations (resp.
permutations). Recall that $J_0=(0,\infty)\subset \mathbb R$. Note
that the action of $\mathbb Z_n$ (resp. $\mathbb P_n$) on $\D$ is
faithful and proper. Let $v\in \D(J_0)$ be a unitary such that
$\beta_g(v)=e^{\frac{2\pi i}{n}} v$ (such $v$ exists by P. 48 of
\cite{ILP}) where $g$ is the generator of the cyclic group $\mathbb
Z_n$ and $\beta_g$ stands for the action of $g$ on $\D$. Note that
$\sigma:=\Ad_v$ is a DHR representation of $\B$ localized on $J_0$.
Let $\gamma: \D(J_0)\rightarrow \B(J_0)$ be the canonical
endomorphism from $ \D(J_0)$ to $\B(J_0)$ and let $\gamma_{\B}:=
\gamma\res \B(J_0)$. Note $[\gamma]=[1]+[g]+...+[g^{n-1}]$ as
sectors of $\D(J_0)$ and
$[\gamma_{\B}]=[1]+[\sigma]+...+[\sigma^{n-1}]$ as sectors of
$\B(J_0)$. Here $[g^i]$ denotes the sector of $\D(J_0)$ which is the
automorphism induced by $g^i.$ All the sectors considered in the
rest of  this paper will be sectors of  $ \D(J_0)$ or  $ \B(J_0)$ as
should be clear from their definitions. All DHR representations will
be assumed to be localized on $J_0$ and have finite statistical
dimensions unless noted otherwise. For simplicity of notations, for
a DHR representation $\sigma_0$ of $\D$ or $\B$ localized on $J_0$,
we will use the same notation $\sigma_0$ to denote its restriction
to $ \D(J_0)$ or  $ \B(J_0)$ and we will make no distinction between
local and global intertwiners for DHR representations localized on
$J_0$ since they are the same by the strong additivity of $\D$ and
$\B$. The following is Lemma 8.3 of \cite{LX}:
\begin{lemma}\label{grading}
Let $\mu$ be an irreducible DHR representation of $\B$. Let
$i$ be any integer. Then:\par
(1) $G(\mu,\sigma^i):=\e (\mu,\sigma^i) \e (\sigma^i,\mu) \in {\mathbb C},
$
$G(\mu,\sigma)^i=G(\mu,\sigma^i) $. Moreover $G(\mu,\sigma)^n=1$;\par
(2) If $\mu_1\prec \mu_2\mu_3$ with
$\mu_1,\mu_2,\mu_3$ irreducible, then $G(\mu_1,\sigma^i)=G(\mu_2, \sigma^i)
G(\mu_3,\sigma^i)$;\par
(3)   $\mu$ is untwisted if and only if  $G(\mu,\sigma)=1;$ \par
(4) $G(\bar \mu, \sigma^i)=\bar G(\mu, \sigma^i).$
\end{lemma}

\subsection{ One cycle case}
First we recall the construction of solitons  for permutation
orbifolds in \S6 of \cite{LX}. Let
$h:S^1\setminus\{-1\}\simeq\mathbb R\to S^1$ be a smooth,
orientation preserving, injective map which is smooth also at
$\pm\infty$, namely the left and right limits $\lim_{z\to
-1^{\pm}}\frac{{\rm d}^n h}{{\rm d}z^n}$ exist for all $n$.

The range $h(S^1\setminus\{-1\})$ is either $S^1$ minus a
point or a (proper) interval of $S^1$.

With $I\in\I$, $-1\notin I$, we set
\[
\Phi_{h,I}\equiv \Ad U(k)\ ,
\]
where $k\in\Diff(S^1)$ and $k(z)=h(z)$ for all $z\in I$ and $U$ is
the projective unitary representation of $\Diff(S^1)$ associated
with $\A$.
Then $\Phi_{h,I}$ does not depend on the choice of
$k\in\Diff(S^1)$ and
\[
\Phi_{h}:I\mapsto \Phi_{h,I}
\]
is a well defined soliton of $\A_0\equiv\A\restriction\mathbb R$.

Clearly $\Phi_h(\A_0(\mathbb R))''=\A(h(S^1\setminus\{-1\}))''$,
thus $\Phi_h$ is irreducible if the range of $h$ is
dense, otherwise it is a type III factor representation.  It is easy
to see that, in the last case, $\Phi_h$ does not depend on $h$
up to unitary equivalence.

Let now $f:S^1\to S^1$ be the degree $n$ map $f(z)\equiv z^n$.
There are $n$ right inverses
$h_i$, $i=0,1,\dots n-1$, for $f$ ($n$-roots); namely there are $n$ injective
smooth maps $h_i:S^1\setminus\{-1\}\to S^1$ such that $f(h_i(z))=z$,
$z\in S^1\setminus\{-1\}$. The $h_i$'s are smooth also at $\pm\infty$.

Note that the ranges $h_i(S^1\setminus\{-1\})$ are $n$ pairwise
disjoint intervals of $S^1$, thus we may fix the labels of the
$h_i$'s so that these intervals are counterclockwise ordered, namely
we have $h_0(1)<h_1(1)<\dots<h_{n-1}(1)<h_0(1)$, and we choose $h_j=
e^{\frac{2\pi ij}{n}}h_0, 0\leq j\leq n-1.$ When no confusion
arises, we will write $h_0$ simply as $z^{\frac{1}{n}}$ and
$\Phi_{h_j, I}(x)= R_{\frac{2\pi j}{n}}R_{z^{\frac{1}{n}}}(x).$
\par
For any interval $I$ of $\mathbb R$, we set
\begin{equation}\label{ts}
\pi_{1,\{0,1...n-1 \},I}\equiv \chi_I\cdot
(\Phi_{h_{0},I}\otimes\Phi_{h_1,I}\otimes\cdots\otimes\Phi_{h_{n-1},I})\ ,
\end{equation}
where $\chi_I$ is
the natural isomorphism from $\A(I_0)\otimes\cdots\otimes\A(I_{n-1})$ to
$\A(I_0)\vee\cdots\vee\A(I_{n-1})$ given by the split property,
with $I_k\equiv h_k(I)$.
Clearly $\pi_{1,\{0,1...n-1 \}} $ is a soliton of
$\D_0\equiv\A_0\otimes\A_0\otimes\cdots\otimes\A_0$ ($n$-fold
tensor product).
Let $p\in \mathbb P_n$. We set
\begin{equation}\label{1cycle1}
\pi_{1,\{p(0),p(1),...,p(n-1) \}}=\pi_{1,\{0,1...,n-1\}} \cdot\b_{p^{-1}}
\end{equation}
where $\b$ is the natural action of $\mathbb P_n$ on $\D$, and
$\pi_{1,\{0,1...,n-1\}} $ is as in (\ref{ts}).
Let $\lambda$ be a DHR representation of $\A$. Given an
interval $I\subset S^1\setminus\{-1\}$,  we set
\begin{definition}\label{1cycle}
\[
\pi_{\lambda,\{ p(0),p(1),...,p(n-1)\}, I}(x)=
\pi_{\lambda,J}(\pi_{1, \{ p(0),p(1),...,p(n-1)\}, I}
(x))\
,\quad x\in\D(I)\ ,
\]
where $ \pi_{1, \{ p(0),p(1),...,p(n-1)\}, I} $ is defined as in
(\ref{1cycle1}),
and $J$ is any interval which contains $I_0\cup I_1\cup...\cup I_{n-1}$.
Denote the corresponding soliton by $\pi_{\lambda,\{ p(0),p(1),...,p(n-1)\}}.$
When $p$ is the identity element in $\mathbb P_n$, we will denote the
corresponding soliton by $\pi_{\lambda, n}$.
\end{definition}
\subsection{Some properties of S matrix for general orbifolds}
Let $\A$ be a completely rational conformal net and let $\Gamma$ be
a finite group acting properly on $\A$. By Th. \ref{orb} $\A^\Gamma$
has only finitely many irreducible representations. We use
$\dot\lambda$ (resp. $\mu$) to label representations of $\A^\Gamma$
(resp. $\A$). We will denote the corresponding genus $0$ modular
matrices by $\dot S, \dot T.$  Denote by $\dot\lambda$ (resp. $\mu$)
the irreducible covariant representations of ${\A}^\Gamma$ (resp.
${\A}$) with finite index. Recall that $b_{\mu\dot\lambda}\in
{\mathbb Z}$ denote the multiplicity of representation $\dot\lambda$
which appears in the restriction of representation $\mu$ when
restricting from $ {\A}$ to ${\A}^\Gamma$. $ b_{\mu\dot\lambda} $ is
also known as the branching rules. We have:
\begin{lemma}\label{Smatrix1}
(1) If $\tau$ is an automorphism (i.e., $d(\tau)=1$) then
$S_{\tau(\lambda) \mu} = G_1(\tau,\mu)^* S_{\lambda\mu}$ where
$\tau(\lambda):= \tau\lambda, G_1(\tau,\mu)=\epsilon(\tau,\mu)
\epsilon(\mu,\tau) ;$\par (2) For any $h\in \Gamma$, let
$h(\lambda)$ be the DHR representation $\lambda\cdot\Ad_{h^{-1}}$.
Then $ S_{\lambda\mu}=  S_{h(\lambda)h(\mu)}$;\par (3) If
$[\alpha_{\dot \lambda}] = [\mu \alpha_{\dot \delta}]$, then for any
$\dot\lambda_1, \mu_1$ with $b_{\dot\lambda_1\mu_1}\neq 0$ we have $
\frac{S_{\dot \lambda \dot\lambda_1}}{S_{\dot 1\dot\lambda_1}} =
\frac{S_{\mu\mu_1}}{S_{ 1\mu_1}} \frac{S_{\dot \delta
\dot\lambda_1}}{S_{\dot 1\dot\lambda_1}} $;\par (4)  We can choose
$c_0(\A^\Gamma)$ so that  $c_0(\A^\Gamma)=c_0(\A)$ (cf. Definition
\ref{c0}).

\end{lemma}
\proof (1), (2) follows from Lemma 9.1 of \cite{KLX}. (3) follows
from the proof on Page 182 of \cite{X3} or Lemma 6.4 of \cite{BE}.
\endproof

\subsection{Fusions of solitons in cyclic orbifolds}
Let $\B\subset \D$ be as in section \ref{cyc}. We note that Th. 8.4
of \cite{KLX} gives a list of all irreducible representations of
$\B.$
\begin{remark}\label{conv}
For simplicity we will label the representation
$\pi_{\lambda,g^j,i}$ ($g=(01...n-1)$) by $(\lambda,g^j,i),$ and
when $i=0$ (resp. $j=0$) which stands for the trivial representation
we will denote the corresponding representation simply as $(\lambda,
g^j)$ (resp. $(\lambda,1)$). When $1$ is used to denote the
representation of a net, it will always be the vacuum
representation.
\end{remark}

\begin{lemma}\label{N}
(1) $G(\sigma,(\mu,g))= e^{\frac{2\pi i}{n}}$ ;\par (2) There exists
an automorphism $\tau_{n,\A}, [\tau_{n,\A}^2]=[1]$ such that
$$ S_{(\lambda,1),(\mu,g)}
=\frac{1}{n}S_{\lambda, \tau_{n,\A}(\mu)} .$$ For simplicity we will
denote $\tau_{n,\A}$ simply as $\tau_n$ when the underlying net $\A$
is clear.
\end{lemma}
\proof (1) follows from remark 4.18 in \cite{DX}, and (2) follows
from Lemma 9.3 of \cite{KLX}.
\endproof
\begin{remark}\label{c02}
By (4) of Lemma \ref{Smatrix1}, we can choose $c_0(\B)$ so that
$c_0(\B)= n c_0(\A).$ We will make such choice in the rest of this
paper.
\end{remark}
\section{Squares of conformal nets}
\begin{definition}\label{square}
Let $\A_i, 1\leq i\leq 4$ be four M\"{o}bius covariant nets such
that $A_3\subset A_2, A_2\subset A_1, A_3\subset A_4$ and
$A_4\subset A_1$ are subnets. Then the square
\(
\begin{array}{ccc}
A_2 & \subset & A_1 \\
  \bigcup &  & \bigcup\\
  A_3 & \subset & A_4\\
\end{array}
\) is called a square of M\"{o}bius covariant nets.
\end{definition}
Let $N=nk, g=(123...N), \D:= \A\otimes\A\otimes...\A.$ Then $g^n$ is
$n$ product of $k$ cycles $g_1,...,g_n,$ with $g_{i+1}=(i (i+n)
(i+2n) ... (i+(k-1)n)), 0\leq i \leq n-1. $ The following square of
conformal nets play an important role in this paper:\par
 \(
\begin{array}{ccc}
 \D^{\l g\r}& \subset & \D^{\l g^n\r} \\
  \bigcup &  & \bigcup\\
  \B_1:=\D^{\l g,g_1,...,g_n\r} & \subset & \B_2:=\D^{\l g_1,...,g_n\r} \\
\end{array}
\)
\begin{proposition}\label{id}
(1) We identify $\D^{\l g_1,...,g_n\r}$ with $n$ tensor products of
a $k$-th order cyclic  permutation orbifold
$(\A\otimes...\otimes\A)^{\l h_1\r}$ in a natural way. Then $\D^{\l
g,g_1,...,g_n\r} \subset \D^{\l g_1,...,g_n\r}$ is a cyclic
permutation of order $n$; Denote by $h_2$ the cyclic permutation on
$\D^{\l g_1,...,g_n\r}$ which comes from permutation
$(01...(n-1))(n(n+1)...(n+n-1))...((k-1)n...((k-1)n+n-1))$ of $\D;$
\par (2)$ \a_{((\lambda,h_1),h_2)}^{\B_1\rightarrow \D^{\l g\r}}
\succ (\lambda,g);$\par
(3)$\a_{((\lambda,h_1^i),1)}^{\B_1\rightarrow \D^{\l
g\r}}=(\lambda,g^{ni})$;\par (4) When $(k,n)=1,$
$\a_{((\lambda,1),h_2^k)}^{\B_1\rightarrow \D^{\l
g\r}}\succ(\lambda,g^{k}).$
\end{proposition}
\proof (1) follows from definition. As for (2), we first  show that
$(\lambda,g)$ and $((\lambda,h_1),h_2)$ come from the restriction of
the same soliton of $\D.$ This can be seen from  definition
\ref{1cycle} as follows: $(\lambda,g)$ comes from a soliton of $\D$
defined by:
$$x_0\otimes x_1\otimes...\otimes x_N \in \A(I)\otimes...\A(I)\rightarrow \pi_\lambda (
R_{z^{\frac{1}{N}}}(x_0) \vee
R_{\frac{2\pi}{N}}R_{z^{\frac{1}{N}}}(x_1)\vee...R_{\frac{2\pi(N-1)}{N}}R_{z^{\frac{1}{N}}}(x_{N-1}))
$$

Let $y_i=x_i\otimes x_{n+i}\otimes...\otimes x_{n(k-1)+i}, 0\leq
i\leq n-1.$ Then $((\lambda,h_1),h_2)$ comes from a soliton of $\D$
defined by
\begin{align*}
y_0\otimes y_1\otimes...\otimes y_{n-1}\rightarrow
&\pi_{\lambda,h_1} (R_{z^{\frac{1}{n}}}(y_0) \vee
R_{\frac{2\pi}{n}}R_{z^{\frac{1}{n}}}(y_1)\vee...R_{\frac{2\pi(n-1)}{n}}R_{z^{\frac{1}{n}}}(y_{n-1})
\\
&=\pi_\lambda(R_{z^{\frac{1}{N}}}(x_0) \vee
R_{\frac{2\pi}{N}}R_{z^{\frac{1}{N}}}(x_1)\vee...R_{\frac{2\pi(N-1)}{N}}R_{z^{\frac{1}{N}}}(x_{N-1}))
\end{align*}

where we have used
$$R_{z^{\frac{1}{k}}}R_{\frac{2\pi i}{n}} (x)= R_{\frac{2\pi i}{kn}}
R_{z^{\frac{1}{k}}}(x),
R_{z^{\frac{1}{k}}}R_{z^{\frac{1}{n}}}(x)=R_{z^{\frac{1}{nk}}}(x) ,
\forall x \in \A(I)$$ Now by Th. 4.8 of \cite{KLX}, $(\lambda,g)$ is
the component of the above soliton where $g$ acts trivially, and
$((\lambda,h_1),h_2)$ is the component of the same soliton where $\l
g,g_1,...,g_n\r$ acts trivially. It follow that the restriction of
$(\lambda,g)$ to $\B_1$ contains $((\lambda,h_1),h_2)$, and (2) is
proved.\par

To prove (3), we first show that
$\a_{((\lambda,h_1^i),1)}^{\B_1\rightarrow \D^{\l g\r}}\succ
(\lambda,g^{ni}). $  As in (2) it is sufficient to show that
$((\lambda,h_1^i),1),(\lambda,g^{ni})$ come from restrictions of the
same soliton of $\D,$ and as in (2) this follows by definition. By
using the index formula in Th. 4.5 and (3) of Prop. 7.4 of
\cite{KLX},  we have $d(((\lambda,h_1^i),1))=d((\lambda,g^{ni})),$
and (3) is proved. (4) is proved in a similar way as in (2):  we
check that $\a_{((\lambda,1),h_2^k)}^{\B_1\rightarrow \D^{\l
g\r}}\succ(\lambda,g^{k})$ by showing that $((\lambda,1),h_2^k),
(\lambda,g^{k})$ come from the same soliton of $\D.$  By definition
\ref{1cycle} $((\lambda,1),h_2^k)$ comes from a soliton of $\D$
defined by
\begin{align*}
x_0\otimes x_1\otimes...\otimes x_N \in \A(I)\otimes...\A(I)
&\rightarrow \pi_\lambda ( R_{z^{\frac{1}{n}}}(x_0) \vee
R_{\frac{2\pi}{n}}R_{z^{\frac{1}{n}}}(x_{-k})\vee... \\
&R_{\frac{2\pi(n-1)}{n}}R_{z^{\frac{1}{n}}}(x_{-k(n-1)}))
\otimes  \\
&\pi_\lambda ( R_{z^{\frac{1}{n}}}(x_{-1}) \vee
R_{\frac{2\pi}{n}}R_{z^{\frac{1}{n}}}(x_{-1-k})\vee...
\\
&R_{\frac{2\pi(n-1)}{n}}R_{z^{\frac{1}{n}}}(x_{-1-k(n-1)})
)\otimes\\
& ...\otimes \pi_\lambda ( R_{z^{\frac{1}{n}}}(x_{-k+1}) \vee \\
&R_{\frac{2\pi}{n}}R_{z^{\frac{1}{n}}}(x_{-k+1-k})\vee...\\
&R_{\frac{2\pi(n-1)}{n}}R_{z^{\frac{1}{n}}}(x_{-k+1-k(n-1)}))
\end{align*}
where indices are defined modulo $N.$ Let $y_{ki}=x_{ki}\otimes
x_{n+ki}\otimes...\otimes x_{n(k-1)+ki}, 0\leq i\leq n-1.$   Since
$(n,k)=1, $ $(\lambda,g^{k})$ comes from a soliton of $\D^{\l
g_1,...g_n\r}$ defined by
\begin{align*}
\pi_{\lambda, (0,k,2k,...,k(n-1))}(y_0\otimes y_1\otimes...\otimes
y_{n-1}) &= \pi_{\lambda, (0,1,2,...,(n-1))}(y_0\otimes
y_{-k}\otimes...\otimes y_{-k(n-1)}) \\
&=\pi_\lambda(R_{z^{\frac{1}{n}}}(y_0)R_{\frac{2\pi}{n}}R_{z^{\frac{1}{n}}}(y_{-k})
\vee... R_{\frac{2\pi (n-1)}{n}}R_{z^{\frac{1}{n}}}(y_{-(n-1)k}))
\end{align*}
where the indexes are defined modulo $n.$  Then the soliton of
$\D^{\l g_1,...g_n\r}$ above comes from restriction of soliton of
$\D$ defined by
\begin{align*}
x_0\otimes x_1\otimes ...\otimes x_{N-1}\rightarrow
& \pi_\lambda (
R_{z^{\frac{1}{n}}}(x_0) \vee
R_{\frac{2\pi}{n}}R_{z^{\frac{1}{n}}}(x_{-k})\vee...R_{\frac{2\pi(n-1)}{n}}R_{z^{\frac{1}{n}}}(x_{-k(n-1)}))
\otimes  \\
& \pi_\lambda ( R_{z^{\frac{1}{n}}}(x_{-n}) \vee
R_{\frac{2\pi}{n}}R_{z^{\frac{1}{n}}}(x_{{-n}-k})\vee
...R_{\frac{2\pi(n-1)}{n}}
R_{z^{\frac{1}{n}}}(x_{{-n}-k(n-1)})) \otimes  \\
& ...\otimes \pi_\lambda ( R_{z^{\frac{1}{n}}}(x_{n(-k+1)}) \vee \\
& R_{\frac{2\pi}{n}}R_{z^{\frac{1}{n}}}(x_{n(-k+1)-k})\vee ...
R_{\frac{2\pi(n-1)}{n}}R_{z^{\frac{1}{n}}}(x_{n(-k+1)-k(n-1)}))
\end{align*}
which up to unitary equivalence (the unitary is a  permutation of
the tensor factors in the Hilbert space) is the same as the soliton
defined at the beginning of the proof of (4). Thus we have shown
that
$$
\a_{((\lambda,1),h_2^k)}^{\B_1\rightarrow \D^{\l
g\r}}\succ(\lambda,g^{k}).
$$
\subsection{Constraints on certain automorphisms}
For simplicity of notations  we define $\tau_{k,n}:=\tau_{n,
(\A\otimes\A\otimes...\otimes \A)^{\mathbb{Z}_k}}.$
\begin{proposition}\label{t1}
(1)$\tau_{n,\A\otimes\A...\otimes \A}= \tau_{n,\A}\otimes
\tau_{n,\A}...\otimes \tau_{n,\A}$ ($k$ tensors);\par
(2)$\tau_{k,n}=
(\tau_n,1,j_{k,n})$ with $k|2j_{k,n}.$
\end{proposition}
\proof Ad (1): Consider inclusions of sunbets $\B_2\subset \D^{\l
g^n\r}\subset \D.$  Note that by definition
$$
\alpha^{\B_2\rightarrow \D^{\l g^n\r}}_{(\lambda_1,g_1)\otimes
(\lambda_2,g_2)\otimes...\otimes
(\lambda_n,g_n)}=(\lambda_1\otimes\lambda_2...\otimes\lambda_n,g^n)
$$
By Lemma \ref{Smatrix1} we have
$$
\frac{S_{(\lambda_1,g_1)\otimes (\lambda_2,g_2)\otimes...\otimes
(\lambda_n,g_n),(\mu_1,1)\otimes(\mu_2,1)\otimes...\otimes(\mu_n,1)}}{S_{1\otimes
...\otimes 1, (\mu_1,1)\otimes...\otimes
(\mu_n,1)}}=\frac{S_{(\lambda_1\otimes...\otimes\lambda_n,g^n),(\mu_1\otimes...\mu_n,1)}}{S_{1\otimes
...\otimes 1, \mu_1\otimes...\otimes \mu_n}}
$$
By Lemma  \ref{N} it follows that
$$
S_{(\tau_{k,\A\otimes...\otimes \A} (\lambda_1
\otimes...\otimes\lambda_n)),\mu_1\otimes...\otimes\mu_n}=S_{(\tau_{k,\A}\lambda_1
\otimes...\otimes\tau_{k,\A}\lambda_n), \mu_1\otimes...\otimes\mu_n}
$$
By unitarity of $S$ matrix, and by replacing $k$ with $n$, (1) is
proved. \par Ad (2): it is sufficient to show that
$\alpha_{\tau_{k,n}}=\tau_n,$ where the induction is with respect to
the $k$-th cyclic permutation orbifold $(\A\otimes \A...\otimes
\A)^{\mathbb{Z}_k}$ and $\A\otimes \A...\otimes \A$ ($k$ tensors).
First we note that since $d(\tau_{k,n})=1,$ and any twisted
representation of $(\A\otimes \A...\otimes \A)^{\mathbb{Z}_k}$ has
index greater or equal to  $k^2$ by Th. 4.5 and Prop. 7.4 of
\cite{KLX}, it follows that $\alpha_{\tau_{k,n}}=\beta$ is a DHR
representation of $\A\otimes \A...\otimes \A$ ($k$ tensors).
 Consider the
square of nets: \(
\begin{array}{ccc}
  \D^{\l h\r}& \subset & \D \\
  \bigcup &  & \bigcup\\
  \B_1=\D^{\l g,g_1,...,g_n\r} & \subset & \B_2:=\D^{\l g^n\r} \\

\end{array}
\) where
$$h=(012...n-1)(n (n+1)... (n+n-1))...(((k-1)n) ((k-1)n+1)...
((k-1)n +n-1))$$ By definition $\alpha^{\B_1\rightarrow \D^{\l
h\r}}_{((\lambda,1),h_2)}= (\lambda,h)$ where by a slight abuse of
notations we use $\lambda$ to denote an irreducible DHR
representation of $\A\otimes...\otimes\A$ ($k$ tensors). By using
(1), Lemma \ref{N} and Lemma \ref{Smatrix1} we have $S_{\beta
\lambda,\mu}= S_{\tau_n \lambda,\mu}$ for all $\lambda,\mu.$ By
unitarity of $S$ matrix and the fact that $[\tau_{k,n}^2]=[1]$ (2)
is proved.
\endproof

\begin{proposition}\label{t2}
(1)$\tau_{k,n}(\lambda,h_1)=(\tau_N\tau_k\lambda,h_1,j)$ for some
$0\leq j\leq k-1;$
\par(2)$\tau_N$ is the vacuum if $N$ is odd, and
$\tau_N=\tau_2$ is $N$ is even;\par (3) $\tau_{k,n}$ is the vacuum
if $n$ is odd, and $j_{k,n}$ as in Prop. \ref{t1} is $0$ modulo $k$
if $k$ is odd.
\end{proposition}
\proof Ad (1): By Prop. \ref{t1} we can assume that
$\tau_{k,n}(\lambda,h_1)=(\mu,h_1,j).$ By Prop. \ref{id} and Lemma
\ref{Smatrix1} we have $S_{(\lambda,g)(\lambda_1,1)}=
S_{\tau_{k,n}(\lambda,h_1),(\lambda_1,1)} ,$ hence $ S_{\tau_k
\mu,\lambda_1}= S_{\tau_N\lambda,\lambda_1}$ by  (2) of Lemma
\ref{N}.  By unitarity of $S$ matrix, (1) is proved.\par Ad (2): By
(1) we have
$\alpha_{\tau_{k,n}(\lambda,h_1)}=\alpha_{(\tau_N\tau_k\lambda,h_1,j)},$
where the induction is with respect to the $k$-th order cyclic
permutation orbifold of $\A\otimes \A\otimes...\otimes\A$ (k
tensors). Note that $\alpha_{\tau_{k,n}}= (\tau_n,...,\tau_n)$ by
Prop. \ref{t1}. It follows by Th. 8.6 of \cite{LX} that
$\tau_n^k=\tau_N\tau_k.$ Hence if $k$ is even,
$\tau_N=\tau_k=\tau_2,$ and if $k$ is odd, $\tau_{nk}= \tau_n
\tau_k.$ Choose $n$ even we have $\tau_k$ is the vacuum when $k$ is
odd.\par Ad (3): (3) follows from (2) and (2) of Prop. \ref{t1}
\endproof
\begin{remark}
We can actually show that $\zeta_k^{j_{k,n}}= \zeta_2^{j_{2,2}}$
when $k,n$ are even integers, but this fact will not be used in the
paper.
\end{remark}
\begin{theorem}\label{94}
$$
[\pi_{1,\{0,1,...,n-1\}}^n]=\bigoplus_{\lambda_1,...,\lambda_n}
M_{\lambda_1,...,\lambda_n} [(\lambda_1,...,\lambda_n)]
$$ where
$M_{\lambda_1,...,\lambda_n}:=\sum_{\lambda} S_{1,\lambda}^{2-2g}
\prod_{1\leq i\leq n} \frac{S_{\lambda_i,\lambda}}{ S_{1,\lambda}} $
with $g=\frac{(n-1)(n-2)}{2},$ and $\pi_{1,\{0,1,...,n-1\}}$ is the
soliton defined in equation (\ref{1cycle1}).
\end{theorem}
\proof This  is proved in Prop. 9.4 of \cite{KLX} under the
assumption $\tau_n^n=1.$ The assumption follows by Prop. \ref{t2}.
\endproof
We note that the above Theorem was conjectured on Page 759 of
\cite{KLX} as a consequence of another conjecture on Page 758 of
\cite{KLX} which states that $\tau_N$ is the vacuum for all $N.$ By
Prop. \ref{t2} it is now enough to prove that $\tau_2$ is the
vacuum.
\subsection{Properties of certain matrices}
In this section we define and examine properties of certain matrices
motivated by P.Bantay's  $\Lambda$ matrices in \cite{B1} and
\cite{B2} which we recalled here for comparison. See \cite{B1} and
\cite{B2} for more details. Suppose that a representation of the
modular group $\Gamma(1)$ has been given.  Let \( r=\frac{k}{n} \)
be a rational number in reduced form, i.e. with \( n>0 \) and \( k
\) and \( n \) coprime. Choose integers \( x \) and \( y \) such
that \( kx-ny=1 \), and define \( r^{*}=\frac{x}{n} \). Then \(
m=\left(
\begin{array}{cc}
k & y\\
n & x
\end{array}\right)  \) belongs to \( \Gamma (1) \), and we define the matrix \( \Lambda \left( r\right)  \)
via \[ \Lambda \left( r\right) _{p,q}=\omega _{p}^{-r}M_{p,q}\omega
_{q}^{-r^{*}}\]
One should fix some definite branch of the logarithm
to make the above definition meaningful, but different choices lead
to equivalent results. See the remark after Lemma \ref{Sprop} for
our choice for genus $0$ modular matrices.

It is a simple matter to show that \( \Lambda (r) \) is well
defined, i.e. does not depend on the actual choice of \( x \) and \(
y \), and  \( \Lambda (r) \) is periodic in \( r \) with period 1,
i.e.
\[ \Lambda (r+1)=\Lambda (r)\]

For \( r=0 \) we just get back the \( S \) matrix \[ \Lambda (0)=S\]
and for a positive integer \( n \) we have \begin{equation}
\label{Lambda1n} \Lambda \left( \frac{1}{n}\right)
=T^{-\frac{1}{n}}S^{-1}T^{-n}ST^{-\frac{1}{n}}
\end{equation}

Finally, we have \[ \Lambda \left( r^{*}\right) _{p}^{q}=\Lambda
(r)_{q}^{p}\]  \[ \Lambda \left( -r\right)
_{p}^{q}=\overline{\Lambda \left( r\right) _{\overline{p}}^{q}}\]

and the functional equation
\begin{equation}\label{fe}
\Lambda (\frac{-1}{r})= T^{\frac{1}{r}} S T^r \Lambda(r) T^{\hat{r}}
\end{equation}
where $\hat{r}= \frac{1}{kn}.$

\begin{definition}\label{hatL}
When $(i,N)=1,$ we define
$$\widehat{\Lambda}_{\lambda_1,\lambda_2}(\frac{i}{N})= N
S_{(\lambda_1,g),(\lambda_2,g^i)},
\widehat{\Lambda}_{\lambda_1,\lambda_2}(r)= N
S_{\lambda_1,\lambda_2}
$$ where $r$ is any integer.
\end{definition}
We note that it follows from the definition that
$$\widehat{\Lambda}_{\lambda_1,\lambda_2}(\frac{i}{N})=
\widehat{\Lambda}_{\lambda_1,\lambda_2}(\frac{i}{N}+1).
$$
\begin{proposition}\label{g2}
(1)$S_{(\lambda_1,g^{i_1},j_1),(\lambda_2,g^{i_2},j_2)}=
\zeta_N^{-i_1j_2-i_2j_1}S_{(\lambda_1,g^{i_1}),(\lambda_2,g^{i_2})};$\par
(2) If $(i_1,N)=1$ and $\hat{i_1} i_1\equiv 1 \mod  N$, then
$S_{(\lambda_1,g^{i_1}),(\lambda_2,g^{i_2})}=
\frac{1}{N}\widehat{\Lambda}_{\lambda_1,\lambda_2}(\frac{i_2\hat{i_1}}{N})
;$\par (3)
$\widehat{\Lambda}_{\lambda_1,\lambda_2}(r)=\widehat{\Lambda}_{\lambda_2,\lambda_1}(r^*);
$ \par (4) $\widehat{\Lambda}_{\lambda_1,\lambda_2}(1-r)=
\overline{\widehat{\Lambda}_{\bar\lambda_1,\lambda_2}}(r); $\par
\end{proposition}
\proof (1) follows from Lemma \ref{Smatrix1} and \ref{grading}. For
(2), let $h\in S_N$ so that $h g^{i_1}h^{-1}=g.$ Then $hgh^{-1}=
g^{\widehat i_1}.$ By Lemma  \ref{Smatrix1} and definition
\ref{hatL}, we have
$$
S_{(\lambda_1,g^{i_1}),(\lambda_2,g^{i_2})}=S_{(\lambda_1,g),
(\lambda_1,g^{i_1\widehat i_2})}=
\frac{1}{N}\widehat{\Lambda}_{\lambda_1,\lambda_2}(\frac{i_2\widehat{i_1}}{N})
$$
For (3), let $h'\in S_N$ be such that $h'g^ih'^{-1}=g.$ Then
$h'gh'^{-1}= g^{\widehat i},$ and by Lemma \ref{Smatrix1} and the
fact that $S$ is symmetric  we have
$$
S_{(\lambda_1,g^i),(\lambda_2,g)}=S_{(\lambda_1,g),(\lambda_2,g^{\widehat
i})}=S_{(\lambda_2,g^{\widehat i}),(\lambda_1,g)}
$$
This proves (3) by definition. \par For (4), by Prop. 6.1 of
\cite{LX} the conjugate of $(\lambda_1,g^i)$ is $(\bar
\lambda_1,g^{-i}).$ (4) now follows from definition and the property
of $S$ matrix under conjugation.
\endproof
\begin{proposition}\label{ggn}
$S_{(\lambda_1,g),(\lambda_2,g^{ni})}=
S_{((\lambda_1,h_1),h_2),((\lambda_2,h_1^i),1)}$
\end{proposition}
\proof Consider the square of nets \(
\begin{array}{ccc}
\D^{\l g\r}& \subset & \D^{\l g^n\r} \\
  \bigcup &  & \bigcup\\
  \B_1:=\D^{\l g,g_1,...,g_n\r} & \subset & \B_2:=\D^{\l g_1,...,g_n\r} \\
\end{array}\) By Prop. \ref{id} and Lemma \ref{Smatrix1} we have
$$
\frac{S_{((\lambda_2,h_1^i),1),((\lambda_1,h_1),h_2)}}{S_{1,((\lambda_1,h_1),h_2)}}
=\frac{S_{(\lambda_2,g^{ni}),(\lambda_1,g)}}{S_{1,(\lambda_1,g)}}
$$
But
$$
\frac{S_{1,((\lambda_1,h_1),h_2)}}{S^{\B_1}_{11}}=d(((\lambda_1,h_1),h_2))
=d(\lambda_1) k^{n-1} \mu_\A^{\frac{(kn-1)}{2}}
$$
and
$$
\frac{S_{1,(\lambda_1,g)}}{S^{\D^{\l g\r}}_{11}}= d(\lambda_1)
\mu_\A^{\frac{(kn-1)}{2}}
$$
where we have used Th. 4.5 and (3) of Prop. 7.4 in \cite{KLX} in the
calculation above. On the other hand
$$
\frac{1}{S^{\B_1}_{11}}=\sqrt{\mu_{\B_1}}= k^n n \sqrt{\mu_\A},
\frac{1}{S^{\D^{\l g\r}}_{11}}=\sqrt{\mu_\D^{\l g\r}}= N
\sqrt{\mu_\A},
$$
and using these equations we obtain
$$
S_{(\lambda_1,g),(\lambda_2,g^{ni})}=
S_{((\lambda_1,h_1),1),((\lambda_2,1),h_2)}
$$
\endproof

\begin{proposition}\label{kn}
Assume that $(k,n)=1$ and $N=kn.$ Then
$S_{(\lambda_1,g^n),(\lambda_2,g^k)}= \frac{1}{N}
S_{\tau_k\lambda_1,\tau_n\lambda_2}.$
\end{proposition}
\proof Consider the square of nets \(\begin{array}{ccc}
  \D^{\l g\r}& \subset & \D^{\l g^n\r} \\
  \bigcup &  & \bigcup\\
  \B_1:=\D^{\l g,g_1,...,g_k\r} & \subset & \B_2:=\D^{\l g_1,...,g_k\r} \\
\end{array}\)
By Prop. \ref{id} and Lemma \ref{Smatrix1} we have
$$
\frac{S_{((\lambda_1,h_1),1),((\lambda_2,1),h_2)}}{S_{1,((\lambda_2,1),h_2)}}
=\frac{S_{(\lambda_1,g^n),(\lambda_2,g^k)}}{S_{1,(\lambda_2,g^k)}}
$$
But
$$
\frac{S_{1,((\lambda_2,1),h_2)}}{S^{\B_1}_{11}}=d(((\lambda_2,1),h_2))
=d(\lambda_2)^k k^{n-1} \mu_\A^{\frac{k(n-1)}{2}}
$$
and
$$
\frac{S_{1,(\lambda_2,g^k)}}{S^{\D^{\l g\r}}_{11}}= d(\lambda_2)^k
\mu_\A^{\frac{k(n-1)}{2}}
$$
where we have used Th. 4.5 and (3) of Prop. 7.4 in \cite{KLX} in the
calculation above. On the other hand
$$
\frac{1}{S^{\B_1}_{11}}=\sqrt{\mu_{\B_1}}= k^n n \sqrt{\mu_\A},
\frac{1}{S^{\D^{\l g\r}}_{11}}=\sqrt{\mu_\D^{\l g\r}}= N
\sqrt{\mu_\A},
$$
and using these equations we obtain
$$
S_{(\lambda_1,g^n),(\lambda_2,g^k)}=
S_{((\lambda_1,h_1),1),((\lambda_2,1),h_2)}
=\frac{1}{n}S_{(\lambda_1,h_1),(\tau_{k,n}(\lambda_2,1))}
$$
Since $(k,n)=1,$ $\zeta_k^{j_{k,n}}=1$ by Prop. \ref{t2},  and we
have
$$
\frac{1}{n}S_{(\lambda_1,h_1),(\tau_{k,n}(\lambda_2,1))}
=\frac{1}{n}\frac{S_{(\lambda_1,h_1),\tau_{k,n}}}{S_{(\lambda_1,h_1)
,1}} \frac{1}{k}S_{\tau_k \lambda_1,\lambda_2} = \frac{1}{N}
S_{\tau_k\lambda_1,\tau_n\lambda_2}
$$\endproof
To prepare the statement of the main theorem in this section, we
define \begin{definition}\label{g} Let $(k,n)=1.$ Define a function
$g$ with value in $Q$ mod $\mathbb{Z}$ by be the following
equations:
$$
g(0)=1, g(\frac{k}{n})= g(\frac{k}{n}\pm 1),
g(\frac{k}{n})+g(\frac{n}{k})=\frac{-2\pi
i}{24}(c-c_0)(3nk-\frac{n^2+k^2+1}{nk})
$$
where $c$ is the central charge of $\A$ and $c_0$ is as in
definition \ref{c0}.
\end{definition}
Such a function is clearly uniquely determined by the defining
equations. We will give further properties of $g$ in Prop.
\ref{dm2}. Let $\Lambda$ be Bantay's $\Lambda$ matrices as reviewed
at the beginning of this section associated with genus 0 $S,T$
matrices as defined after definition \ref{c0}. Then we have:
\begin{theorem}\label{keyeq}
$$\widehat\Lambda_{\lambda_1,\lambda_2}(r)=\exp(2\pi i
g(r))\Lambda_{\lambda_1,\lambda_2}(r)
$$ where $r\in \mathbb{Q}$ and $g(r)$ is as in definition \ref{g}.
\end{theorem}
\proof The idea is to consider equation $S=TSTST$ for the net
$\D^{\l g\r}$ with the order of $g$ equal to $nk$ and $(n,k)=1.$

Let us compute the $(\lambda_1,g^n),(\lambda_2,g^{k-n})$ entry on
both sides of the equation. Since $(k,n)=1,$ $(k-n,kn)=1.$ Let
$x_1,x_2$ be integers such that $x_1(k-n)+knx_2=1.$

By Prop. \ref{ggn} the left hand side is

$$
S_{(\lambda_1,g^n),(\lambda_2,g^{k-n})} =
S_{(\lambda_1,g^{nx_1}),(\lambda_2,g)}
=S_{((\lambda_1,h_1^{x_1}),1), ((\lambda_2,h_1),h_2)}
$$
By Lemma \ref{N} and Prop. \ref{t2} we have

$$
S_{((\lambda_1,h_1^{x_1}),1), ((\lambda_2,h_1),h_2)} =\frac{1}{n}
S_{(\lambda_1,h_1^{x_1}),\tau_{k,n}(\lambda_2,h_1)} =\frac{1}{N}
\frac{S_{\tau_k(\lambda_1),\tau_n}}{S_{\tau_k(\lambda_1),1}}
\widehat\Lambda_{\lambda_1\lambda_2}(\frac{k-n}{k})
$$

By (14) of \cite{LX} and remark \ref{c02} we have

$$
T_{\lambda_1,g^n} = T_{\lambda_1}^{\frac{n}{k}}\exp({2\pi i
(\frac{(k^2-1)(c-c_0)}{24k})}, T_{\lambda_3,g^k} =
T_{\lambda_3}^{\frac{k}{n}}\exp(2\pi i (\frac{(n^2-1)(c-c_0)}{24n}))
$$
and
$$
T_{\lambda_2,g^{k-n}} = T_{\lambda_2}^{\frac{1}{nk}}\exp(2\pi i
(\frac{(n^2k^2-1)(c-c_0)}{24nk}))
$$
By using the above equations and Prop. \ref{kn} we obtain the
$(\lambda_1,g^n),(\lambda_2,g^{k-n})$ entry on the RHS is given by
$$
\frac{1}{N}\sum_{\lambda_3} \exp(\frac{2\pi i
(c-c_0)(3nk-\frac{n^2+k^2+1}{nk})}{24}) T_{\lambda_1}
S_{\tau_k(\lambda_1),\tau_n(\lambda_3)}T_{\lambda_3}^{\frac{k}{n}}
\widehat\Lambda_{\lambda_3,\lambda_2}(\frac{k-n}{n})
\frac{S_{\tau_k,\tau_n(\lambda_3)}}{S_{1\lambda_3}}
T_{\lambda_2}^{\frac{1}{kn}}
$$
Since $(k,n)=1$, by Prop. \ref{t2} we have
$$
\frac{S_{\tau_k,\lambda_3}}{S_{1,\lambda_3}}=\frac{S_{\tau_k,\tau_n\lambda_3}}{S_{1,\lambda_3}}
$$
when $k$ is odd, and when $k$ is even, $n$ must be odd and the above
equation also holds. When comparing both the LHS and RHS, we see
that the $\tau$ dependence canceled out from both sides and we are
left with the following equation for $\widehat\Lambda:$
$$
\widehat\Lambda_{\lambda_1,\lambda_2}(\frac{k-n}{k})=
\exp(\frac{2\pi i
(c-c_0)(3nk-\frac{n^2+k^2+1}{nk})}{24})T_{\lambda_1}
S_{\lambda_1,\lambda_3}T_{\lambda_3}^{\frac{k}{n}}
\widehat\Lambda_{\lambda_3,\lambda_2}(\frac{k-n}{n})
T_{\lambda_2}^{\frac{1}{kn}}
$$
Comparing with the   equation (\ref{fe}) of Bantay's $\Lambda$
matrices and using Prop. \ref{g2} we conclude that there is a $\mod
\mathbb{Z}$ valued function as defined in definition \ref{g} such
that
$$\widehat\Lambda_{\lambda_1,\lambda_2}(r)=\exp(2\pi i
g(r))\Lambda_{\lambda_1,\lambda_2}(r).
$$
\endproof
Note that the theorem above determined $\widehat{\Lambda}$ matrices
completely, and hence the entries of $S$ matrix as in definition
\ref{hatL}\footnote{With little effort we can in fact determine all
entries of $S$ matrix for the cyclic permutation orbifold using the
methods of this chapter. However Th. \ref{keyeq} is enough for the
purpose of this paper.}. By using Verlinde's formula, one can write
down a series of equations of fusion rules in terms of
$\widehat{\Lambda}$ matrices. Since fusion coefficients are
non-negative integers, these equations describe certain arithmetic
properties of $\widehat{\Lambda}$ matrices, and none of them seems
to trivial for the case of conformal nets associated with $SU(n)$ at
level $k$ where $S,T$ matrices are given (cf. \cite{W}). We refer
the reader to Cor. 9. 9 for such a statement in the  case when
$N=2.$

\section{Arithmetic properties of $S,T$ matrices for a completely rational net}

\subsection{Galois action on \protect\( \widehat{\Lambda} \protect \) matrices}

In this section we'll study the Galois action in the cyclic
permutation orbifold $\D^{\l g\r}$ as in \cite{B1}.  By  Th.
\ref{orb} the Galois action on the genus 0 \( S \)-matrix elements
of $\D^{\l g\r}$  may be described via suitable permutations \( \pi
_{l} \) of the irreducible representations of the orbifold and signs
\( \varepsilon _{l}. \) This will in turn allow us to determine the
Galois action on \( \widehat{\Lambda }\)-matrices as defined in
definition \ref{hatL}.

Let $N$ be  a positive integer , and as in \S3 consider  the cyclic
permutation \(  g=( 1,\ldots ,N ) \). We will use $C(N,\A)$ to
denote the conductor of the permutation orbifold $\D^{\l g\r}.$
Hence $C(1,\A)$ is the conductor of $\A.$ Note that $C(1,\A)$
depends on the choice of $c_0(\A)$ in definition \ref{c0}, and our
choice of $c_0(\D^{\l g \r})$ is as in remark \ref{c02}.\par

Among the irreducible representations of the permutation orbifold
$\D^{\l g\r}$ there is a subset \( \mathcal{J} \) of special
relevance to us. The elements in \( \mathcal{J} \) are labeled by
triples \( \left( \lambda,g^n,k\right) \), where \( \lambda  \) is
an irreducible representation of \( \mathcal{A} \), while \( n \)
and \( k \) are integers mod \( N \). The subset of those \( \left(
\lambda,g^n,k\right) \) where \( n \) is coprime to \( N \) will be
denoted by \( \mathcal{J}_{0} \). It follows from (1) of Lemma
\ref{Smatrix1} that  \( \left( \lambda, g^n,k\right) \in
\mathcal{J}_{0} \) have vanishing \( S \)-matrix elements with the
labels not in \( \mathcal{J} \), while for \( \left(
\mu,g^m,l\right) \in \mathcal{J} \) we have:
\begin{equation}
\label{Smat} S_{\left( \lambda,g^n,k\right),\left( \mu,g^m,l\right)
}=\frac{1}{N}\zeta _{N}^{-(km+ln)}\widehat{\Lambda}
_{\lambda,\mu}\left( \frac{m\widehat{n}}{N}\right)
\end{equation}
where \( \widehat{n} \) denotes the mod \( N \) inverse of \( n \)
and \( \zeta _{N}=\exp \left( \frac{2\pi i}{N}\right)  \).
\begin{lemma}\label{ti}
$\varepsilon_l(\tau_N(\lambda))= \varepsilon_l (\lambda),
\tau_N\pi_l(\lambda) = \pi_l(\tau_N\lambda). $
\end{lemma}
\proof Note that by the property of $\tau_N$ we have
$\frac{S_{\tau_N,\mu}}{S_{1,\mu}}=\pm 1.$ By definition of Galois
actions we have
$$
\sigma_l(S_{\tau_N\lambda,\mu})=\frac{S_{\tau_N,\mu}}{S_{1,\mu}}\sigma_l(S_{\lambda,\mu})
=\frac{S_{\tau_N,\mu}}{S_{1,\mu}}\varepsilon_l(\lambda)
S_{\pi_l(\lambda),\mu} =\varepsilon_l(\tau_N \lambda)
S_{\pi_l(\tau_N\lambda),\mu}
$$
Hence
$$
\varepsilon_l(\lambda) S_{\tau_N\pi_l(\lambda),\mu}=
\varepsilon_l(\tau_N(\lambda)) S_{\pi_l(\tau_N\lambda),\mu}
$$
By unitarity of $S$ matrix the lemma is proved.
\endproof
By using Lemma \ref{ti} and Th. \ref{keyeq}, the proofs of Lemma 1-3
Prop.1, Cor. 1 and Th. 1 of \cite{B1} go through (In the statements
of Lemma 1-3, Prop. 1,  Bantay's $\Lambda$ matrix has to be replaced
by our $\widehat{\Lambda}$ matrix, and the additional assumption on
$l$ is that $l$ is coprime to $C(N,\A)$ where $N$ is the denominator
of a rational number $r$ as given in these statements)

For reader's convenience and to set up notations, we summarize
Lemmas 1-3 and Prop. 1 of \cite{B1} in the following and sketch its
proof.
\begin{lemma}\label{Z}
Assume that $l$ is coprime to the denominator $N$ of $r$ and
$C(N,\A)C(1,\A)$. Then:\par (1)  The set \( \mathcal{J} \) is
invariant under the permutations \( \tilde{\pi }_{l} \), i.e. \(
\tilde{\pi }_{l}\left( \mathcal{J}\right) =\mathcal{J} \). For  \(
\left( \lambda,g^n,k\right) \in \mathcal{J}_{0} \) one has
\begin{equation} \label{pi} \tilde{\pi }_{l}\left(
\lambda,g^n,k\right) =\left( \pi
_{l}(\lambda),g^{ln},\tilde{k}\right)
\end{equation}
for some function \( \tilde{k} \) of \( l,\lambda,n \) and \( k \),
and
\begin{equation} \label{eps} \tilde{\varepsilon }_{l}\left(
\lambda,g^n,k\right) =\varepsilon _{l}(\lambda)
\end{equation}\par
(2) \begin{equation} \label{Lmat2} \sigma _{l}\left(
\widehat{\Lambda} \left( r\right) \right) =\widehat{\Lambda} \left(
lr\right) G_{l}Z_{l}(r^{*})=Z_{l}\left( r\right)
G_{l}^{-1}\widehat{\Lambda} (\widehat{l}r)
\end{equation}
where \( Z_{l}(r) \) is a diagonal matrix whose order divides the
denominator \( N \) of \( r \), and \( \widehat{l} \) is the mod \(
N \) inverse of \( l \), and  \( Z_{l}(0)=\mathbb {I} \) , \(
Z_{l}(r+1)=Z_{l}(r) \);\par (3)
\begin{equation}
\label{zcoc}
G_{l}^{-1}Z_{m}(\widehat{l}r)G_{l}=Z_{lm}(r)Z_{l}^{-m}(r)
\end{equation}
whenever both \( l \) and \( m \) are coprime to the denominator $N$
of \( r \) and $C(N,\A)C(1,\A)$;\par (4)  If \( n \) is coprime to
the denominator of \( r \), then
\begin{equation} \label{zmult} Z_{l}^{n}\left( r\right) =Z_{l}\left(
nr\right)
\end{equation}
\end{lemma}
\proof We give a proof of (1) following the proof of Bantay
indicating necessary changes. First, let's fix \( \left(
\lambda,g^n,k\right) \in \mathcal{J} \). According to
Eq.(\ref{Smat}), we have \[ S_{\left( \lambda,g^n,k\right),\left(
\mu,1\right) }=\frac{1}{N}\zeta _{N}^{-k} \widehat{\Lambda}
_{\lambda,\mu}\left( \frac{\hat{n}}{N}\right) \] and this expression
differs from \( 0 \) for at least one \( \mu \), by the unitarity of
\( \Lambda \)-matrices and Th. \ref{keyeq}. Select such a \( \mu \),
and apply \( \sigma _{l} \) to both sides of the equation. One gets
that \[ \tilde{\varepsilon }_{l}\left( \lambda,g^n,k\right)
S_{\tilde{\pi }_{l}\left( \lambda,g^n,k\right),(\mu,1) }=\sigma
_{l}\left( S_{\left(\lambda,g^n,k\right),(\mu,1) }\right)
\] differs from \( 0 \), but this can only happen if \( \tilde{\pi
}_{l}\left( \lambda,g^n,k\right] \in \mathcal{J} \) because \(
\left( \mu,1\right) \in \mathcal{J}_{0} \).

Next, for \( \left[ \lambda,g^n,k\right] \in \mathcal{J}_{0} \)
consider\begin{equation} \label{S0} S_{\left(
\lambda,g^n,k\right),(\mu,1,m) }=\frac{1}{N}\zeta
_{N}^{-nm}S_{\tau_N\lambda,\mu}
\end{equation}
Applying \( \sigma _{l} \) to both sides of the above equation we
get from Eq.(\ref{sls}) \[ \tilde{\varepsilon
}_{l}(\lambda,g^n,k)S_{\tilde{\pi }_{l}\left( \lambda,g^n,k\right),
(\mu,1,m) }=\frac{1}{N}\zeta _{N}^{-lnm}\varepsilon _{l}(\tau_N
\lambda)S_{\pi _{l}(\tau_N \lambda), \mu}\] But the lhs. equals \[
\tilde{\varepsilon }_{l}(\lambda,g^n,k)\frac{1}{N}\zeta
_{N}^{-\tilde{n}m}S_{\tau_N\tilde{\lambda},\mu}\]
 according to Eq.(\ref{S0}) if \( \tilde{\pi }_{l}\left( \lambda,g^n,k\right)
 =\left[ \tilde{\lambda},g^{\tilde{n}},\tilde{k}\right]  \).
Equating both sides  we arrive at \[
S_{\tau_N\tilde{\lambda},\mu}=\varepsilon
_{l}(\tau_N\lambda)\tilde{\varepsilon }_{l}(\lambda,g^n,k)\zeta
_{N}^{-m(\tilde{n}-ln)}S_{\pi _{l}(\tau_N\lambda),\mu}\] By Lemma
\ref{ti} and the fact that
$$
S_{\tau_N\lambda,\mu}= \frac{S_{\tau_N,\mu}}{S_{1,\mu}}
S_{\lambda,\mu}
$$
we have
\[
S_{\tilde{\lambda},\mu}=\varepsilon _{l}(\lambda)\tilde{\varepsilon
}_{l}(\lambda,n,k)\zeta _{N}^{-m(\tilde{n}-ln)}S_{\pi
_{l}(\lambda),\mu}\]

As the lhs. is independent of \( m \), we must have \[
{\tilde{n}}={ln} \ \mod {N}\] and \[ \tilde{\lambda}=\pi
_{l}(\lambda)\]
 as well as \[
\tilde{\varepsilon }_{l}(\lambda,g^n,k)=\varepsilon _{l}(\lambda)\]
The proof of (2)-(4) is the same  as that of Bantay with his
$\Lambda$ matrices replaced with our $\widehat{\Lambda}$ matrices.
\endproof

 Let us sketch the proof of the following (cf. Th. 1 in
\cite{B1}) theorem, indicating modifications compared to the proof
in \cite{B1}:
\begin{theorem}\label{th1}
Let $\A$ be a completely rational net and let $T$-matrix be defined
as after definition \ref{c0}.  Then for all $l$ coprime to the
conductor $G_l^{-1} T G_l= T^{l^2}.$
\end{theorem}
\proof Let $N$ be the order of $T.$ Then $N$ divides the conductor
by definition. Choose $l$ so that $(l,12 C(1,\A)C(N,\A))=1.$ By Th.
\ref{keyeq} we have $\widehat{\Lambda}(\frac{1}{N})=
\exp(\frac{-2\pi i(c-c_0)(N^2-1)}{12N}) \Lambda (\frac{1}{N}).$
Follow the argument  of Bantay, with $\Lambda$ replaced by
$\widehat{\Lambda}$ we have
$$
\exp(\frac{-2\pi i(c-c_0)(N^2-1)(l^2-1)}{12N}) T^{\frac{-2l^2}{N}} =
G_l^{-1} T^{\frac{-2}{N}} G_l Z_l^2(\frac{l}{N})
$$
Now we use the fact that since $(l,12)=1, 12|l^2-1.$ The rest of the
argument is the same as \cite{B1} and we have $G_l^{-1} T G_l=
T^{l^2}$ for $l$ with the property $(l,12 C(1,\A)C(N,\A))=1.$ Now
for any $l$ coprime to the conductor $C(1,\A),$ by Dirichelet
theorem on arithmetic progressions we can always find integer $p$ so
that $l_1=l+pC(1,\A)$ with the property that $(l,12
C(1,\A)C(N,\A))=1.$ Since $G_{l_1}=G_l, T^{l_1^2}=T^{l^2},$ the
theorem is proved for any $l$ coprime to the conductor.
\endproof
This above theorem has been conjectured in \cite{CG}, where some of
its consequences had been derived. \par
Prop. 2 of \cite{B1} has to
be modified due to phase factor as follows:
\begin{proposition}\label{p2}
Let $r=\frac{n}{N}.$ If $l$ is coprime to $C(N,\A)C(1,\A)N,$ then
$$
G_l^{-1}T^r G_l= T^{l^2 r} Z_l^l(r) \exp(\frac{-2\pi
i(l^2-1)(c-c_0)r}{24})
$$
\end{proposition}
\proof The proof is similar to the proof of Prop. 2 in \cite{B1} and
we indicate modifications when necessary. Write $r=\frac{n}{N}.$ The
idea is to apply Th. \ref{th1} to  $\D^{\l g\r}$ with the order of
$g$ equal to $N.$ The phase factor comes in when we  note that
$$T_{(\lambda,g^n,k)}= \zeta_N^{nk} T_{\lambda}^{\frac{1}{N}}
\exp(\frac{2\pi i(c-c_0)}{24}(N-\frac{1}{N}))$$ Use Th. \ref{th1} we
have
$$
T_\lambda^{\frac{l^2}{N}}\exp(\frac{2\pi
i(c-c_0)(l^2-1)}{24}(N-\frac{1}{N})) =
\zeta_N^{lk_0}T_{\pi_l(\lambda)}^{\frac{1}{N}}
$$
and the rest of the proof is as in
 \cite{B1} .
\endproof

By using Th. \ref{th1}, Prop. 3-6 and Cor. 2 of \cite{B1} follows in
our setting with the same proof (except (4) in the theorem below )as
in \cite{B1} and \cite{CG}. We record these results in the following
theorem:
\begin{theorem}\label{arithmetic1}
Let $\A$ be completely rational net and let $S,T$ be its genus 0
modular matrices as defined after definition \ref{c0}. Then:\par (1)
For \( l \) coprime to the conductor,
\begin{equation} \label{Gl} G_{l}=S^{-1}T^{l}ST^{\widehat{l}}ST^{l}
\end{equation}
where \( \widehat{l} \) denotes the inverse of \( l \) modulo the
conductor;\par (2) The conductor equals the order \( N \) of \( T
\), and \( F=\mathbb{Q}\left[ \zeta _{N}\right] \); \par (3) Let \(
N_{0} \) denote the order of the matrix \( \omega _{0}^{-1}T \),
i.e. the least common multiple of the denominators of the conformal
weights. Then \( N=eN_{0} \), where the integer \( e \) divides \(
12 \). Moreover, the greatest common divisor of \( e \) and \( N_{0}
\) is either 1 or 2;\par (4) \( N_{0} \) times the central charge
$c$ is an even integer; \par (5) There exists a function \( N(r) \)
such that the conductor \( N \) divides \( N(r) \) if the number of
irreducible representations of $\A$ - i.e. the dimension of the
modular representation - is \( r \).
\end{theorem}
\proof Given Th. \ref{th1}, (1),(2), (3), and (5) are proved in the
same way as in \cite{B1}. As for (4), the proof of Cor. 2 in
\cite{B1} shows that $N_0 c_0$ is an even integer. By Lemma 9.7 of
\cite{KLX} $c-c_0\in 4 {\mathbb Z}$ and (4) is proved.
\endproof

\subsection{The kernel of the modular representation}
In this section we consider the modular representation of a
completely rational net $\A$ as defined after Lemma \ref{Sprop}. We
will show that this representation factorizes through a congruence
subgroup. We refer the reader to \cite{CG} for a nice  account of
this and related questions.
 Recall that the kernel \( \mathcal{K}
\) consists of those modular transformations which are represented
by the identity matrix, i.e.
\[ \mathcal{K}=\left\{ m\in \Gamma (1)\, |\, M_{\lambda,\mu}=\delta
_{\lambda,\mu}\right\} \]

\begin{proposition}\label{dm}
If \( m=\left( \begin{array}{cc}
a & b\\
e & d
\end{array}\right) \in \Gamma (1) \) with \( d \) coprime to the
conductor \footnote{We use $e$ instead of more natural $c$ since $c$
has been used to denote the central charge.}. Let $m_e$ be an
integer such that $m_e g(\frac{a}{e})\in \mathbb{Z}.$ Let $l=d
+meC(1,\A)$ be such that $l$ is coprime to $6m_eC(|e|,\A).$ Then
$$
\sigma_d(M) =T^{b}S^{-1}T^{-e}\sigma _{d}(S) \exp({-2\pi
i}(lg(\frac{a}{e})-g(\frac{1}{e})-\frac{(l^2-1)(c-c_0)}{24e}))
$$
\end{proposition}
\begin{proof}
According to Eq.(\ref{Lmat2}) and Th. \ref{keyeq}
$$
\sigma _{l}(M) =\sigma _{l}( T^{a/e}\Lambda( \frac{a}{e}) T^{d/e})
=\sigma _{l}( T^{a/e} \exp(-2\pi i g(\frac{a}{e}))
\widehat{\Lambda}( \frac{a}{e}) T^{d/e})
$$
By our assumption on $l$ we have
\begin{align*}
\sigma _{l}( T^{a/e} \exp(-2\pi i g(\frac{a}{e})) \widehat{\Lambda}(
\frac{a}{e}) T^{d/e})
&= \exp(-2\pi i (l g(\frac{a}{e})-
g(\frac{1}{e})  ))\times \\
& T^{ad/e}\Lambda ( \frac{1}{e}) G_{l}Z_{l}( \frac{d}{e}) T^{d^2/e}
\end{align*}
But

$$
\Lambda( \frac{1}{e}) =T^{-1/e} S^{-1} T^{-c} ST^{-1/e} $$ so
\begin{align*}
\sigma _{l}( M) &=\exp(-2\pi i (l g(\frac{a}{e})-
g(\frac{1}{e})  ))\times \\
& T^b S^{-1} T^{-e} S T^{-1/e} G_lZ_l(d/e) T^{\frac{d^2}{e}}
\end{align*}
From Prop. \ref{p2}
$$
T^{-1/e}G_{l}=G_{l}T^{-l^{2}/e}Z_{l}^l(-1/e) \exp(\frac{2\pi
i(l^2-1)(c-c_0)}{24e})
$$
Putting all this together and using Lemma \ref{Z}  we get the
proposition.
\end{proof}
Next we show that the phase factor in the above proposition is
always $1$:
\begin{proposition}\label{dm2}
Let $l$ be as in Prop. \ref{dm}. Then
$$
\exp({-2\pi
i}(lg(\frac{a}{e})-g(\frac{1}{e})-\frac{(l^2-1)(c-c_0)}{24e})) =1$$
\end{proposition}
\proof Let $4x=(c-c_0).$ By Lemma 9.7 of \cite{KLX} $x$ is an
integer. Let us first prove the proposition for the case $x=2x_2$ is
even. Choose an integer $n$ so that $3n+x_2
>0$ and consider a local net $\E$ which is $3n+x_2$ tensor product
of the local net $\A_{(E_8)_1}.$ This net has $\mu$ index equal to
one by Th. 3.18 of \cite{DX}.  We choose our $c_0(\E)=24n$ in the
definition of $T$ matrix for $\E.$ The corresponding modular
representation is trivial. We denote by $\widehat{\Lambda}_\E$
(resp. $\Lambda_\E$) the matrices as defined in definition
\ref{hatL} (resp. Bantay's $\Lambda$ matrices) associated with $\E$.
Apply Th. \ref{keyeq} to $\E$ we have
$$
\widehat{\Lambda_\E}(a/e)=\exp(2\pi ig(a/e))\Lambda_\E(a/e)
$$
By (2) of Th. \ref{arithmetic1} the conductor of the $e$-th cyclic
permutation orbifold of $\E$ divides $3e.$  By conditions on $l$ we
can apply Prop. \ref{dm} to the net $\E$ to have
$$
\sigma_d(M) =T^{b}S^{-1}T^{-c}\sigma _{d}(S) \exp({-2\pi
i}(lg(\frac{a}{e})-g(\frac{1}{e})-\frac{(l^2-1)(c-c_0)}{24e}))
$$
Since the modular representation for $\E$ is trivial we must have
$$
\exp({-2\pi i}(lg(\frac{a}{e})-g(\frac{1}{e})-\frac{(l^2-1)( c-
c_0)}{24e}))=1
$$
and we have proved proposition for $x$ even.\par If $x=2(x_1+2)-3$
is odd, we define a new  set of $S_1,T_1$ matrix by
$$
T_1=\exp{\frac{-2\pi i x }{6}} T, S_1=\exp{\frac{-2\pi i x}{2}} T
$$
and denote by $\Lambda_1$ the $\Lambda$ matrix of Bantay associated
with $S_1,T_1$. Let $g_1(k/n)$ be defined modulo integers such that
$\Lambda_1 = \exp(2\pi i g_1(k/n)) \Lambda.$ From the definition
\ref{g} one checks easily that modulo integers
$$ g(k/n)= \frac{n}{2} + g_1(k/n)$$
Using the assumption that $l$ is odd it is now sufficient to check
the proposition for $g_1.$ From the defining equation for $g_1$, we
see that $\exp{2\pi i g_1(k/n)}$ is Bantay's $\Lambda$ matrix
associated with one dimensional representation of the modular group
given by $T\rightarrow \exp{\frac{-2\pi i x }{6}}, S\rightarrow
\exp{\frac{-2\pi i x }{2}} .$  This representation is the tensor
product of two one dimensional representations given by
$$
T\rightarrow \exp(\frac{-2\pi i (2x_1+4) }{6}), S\rightarrow 1
$$
and
$$
T\rightarrow \exp{\frac{2\pi i  }{2}}, S\rightarrow \exp{\frac{2\pi
i  }{2}}
$$
and we denote by $g_3,g_2$ the associated $\Lambda$ matrices. Note
that $g_1(a/e)=g_2(a/e)+g_3(a/e).$ The same proof as in the $x$ even
case, with $x_2=x_1+2$, shows that
$$
\exp({-2\pi
i}(lg_3(\frac{a}{e})-g_3(\frac{1}{e})-\frac{(l^2-1)(2x_1+4)}{24e}))
=1
$$
Hence to finish the proof we just have to show that
$$
\exp({-2\pi
i}(lg_2(\frac{a}{e})-g_2(\frac{1}{e})+\frac{(l^2-1)}{2e})) =1$$
Since the associated modular representation is very simple, this can
be checked directly using the following formulas:
$$
\exp(2\pi i g_2(a/e)) = \exp(\frac{-2\pi i(a+d)}{2e})M_2(a,b,e,d)
$$
where $M_2(a,b,e,d)=\pm 1$. When $e$ or $d$ is even,
$M_2(a,b,e,d)=(-1)^d$ or $(-1)^e;$ When $e,d$ are odd,
$M_2(a,b,e,d)=(-1)^{a+d+1}.$
\endproof
By combining the above two propositions we have proved the
following:
\begin{theorem}\label{dm3}
If $d$ is coprime to the conductor, then
$$
\sigma_d(M)= T^b S^{-1} T^{-e} \sigma_d(S).
$$
\end{theorem}
Now Th. 2-4 of \cite{B1} follow exactly in the same way. Let us
record these theorems in the following:
\begin{theorem}\label{cong}
Let $\A$ be a completely rational net and consider the modular
representation as defined after Lemma \ref{Sprop}. Then:\par (1) Let
\( d \) be coprime to the conductor \( N \). Then \( \left(
\begin{array}{cc}
a & b\\
e & d
\end{array}\right) \in \Gamma (1) \) belongs to the kernel \( \mathcal{K} \) if and only if
\begin{equation}
\label{kern1} \sigma _{d}\left( S\right) T^{b}=T^{e}S;
\end{equation} \par (2)
Define \[ \Gamma _{1}\left( N\right) =\left\{ \left(
\begin{array}{cc}
a & b\\
e & d
\end{array}\right) \in \Gamma \left( 1\right) \, |\, {a,d}\equiv {1}\mod \ {N},\, \,
{e}\equiv{0}\mod \ {N}\right\} \] and\[ \Gamma \left( N\right)
=\left\{ \left( \begin{array}{cc}
a & b\\
e & d
\end{array}\right) \in \Gamma _{1}\left( N\right) \, |\, {b}\equiv{0}\mod  {N} \right\} \]
Then\[ \mathcal{K}\cap \Gamma _{1}\left( N\right) =\Gamma \left(
N\right)
\]
 In particular, \( \mathcal{K} \) is a congruence subgroup of level \( N
 \);\par (3)
Define \( SL_{2}\left( N\right) \cong \Gamma (1)/\Gamma (N) \). The
modular representation factorizes through $SL_{2}\left( N\right)$
which we denote by $D$.
For \( l \) coprime to \( N \), define the automorphism \( \tau
_{l}:SL_{2}\left( N\right) \rightarrow SL_{2}\left( N\right)  \) by
\begin{equation} \label{taudef} \tau _{l}\left( \begin{array}{cc}
a & b\\
e & d
\end{array}\right) =\left( \begin{array}{cc}
a & lb\\
\widehat{l}e & d
\end{array}\right)
\end{equation}
where \( \widehat{l} \) is the mod \( N \) inverse of \( l \). Then
\( \label{gal2} \sigma _{l}\circ D=D\circ \tau _{l}. \)
\end{theorem}

\end{document}